\documentclass[10pt]{amsart} 

\usepackage[utf8]{inputenc}
\usepackage[T1]{fontenc}	
\usepackage[french,english]{babel}	

\usepackage{lmodern}			

\usepackage{graphicx}	

\usepackage[a4paper,plainpages=false,colorlinks,linkcolor=bleuFonce,citecolor=rougeFonce,urlcolor=vertFonce,breaklinks]{hyperref}

\usepackage{placeins}
\usepackage{float} 

\usepackage{color}
\definecolor{vertFonce}{rgb}{0,0.5,0}
\definecolor{numLignes}{rgb}{0.17,0.57,0.7}	
\definecolor{gris}{rgb}{0.5,0.5,0.5}
\definecolor{grisFonce}{rgb}{0.2,0.2,0.2}
\definecolor{orange}{rgb}{1,0.65,0.31}		
\definecolor{orangeFonce}{rgb}{1,0.4,0}
\definecolor{bleuFonce}{rgb}{0,0,0.4}
\definecolor{rougeFonce}{rgb}{0.3,0,0}
\definecolor{rougeWord}{rgb}{0.5,0,0}
\definecolor{vertClair}{rgb}{0.8,1,0.8}
\definecolor{rougeClair}{rgb}{1,0.5,0.5}

\usepackage{pict2e}
\usepackage{multido}

\usepackage{amsfonts,amssymb,amsthm,amsmath} 
\usepackage{dsfont}				
\usepackage{mathrsfs}
\usepackage{yfonts}
 
\newtheorem{lem}{Lemma}[section]
\newtheorem{theorem}{Theorem}

\newtheorem{prop}{Proposition}[section]

\newenvironment{demo}[1][]{%
	\begin{proof}[\textbf{Proof of #1}]
	}{%
	\end{proof}
}

\newenvironment{thm}[1][]{%
	\color{rougeWord}\begin{theorem}[#1]
	}{%
	\end{theorem}
}

\newcommand{\step}[1]	{\paragraph{\itshape\bfseries Step #1.}}

\newcommand		{\N}			{\mathbb N}			
\newcommand		{\R}			{\mathbb R}			
\newcommand		{\B}			{\mathscr B}		

\newcommand		{\ssi}			{\Leftrightarrow}

\newcommand		{\oball}		{\mathring B}

\renewcommand	{\d}			{\,\mathrm{d}}	
\DeclareMathOperator{\divg}		{div}

\DeclareMathOperator{\G}		{\Gamma}		
\DeclareMathOperator{\diam}		{diam}
\DeclareMathOperator{\vp}		{vp}

\newcommand		{\weight}[1]	{\langle #1\rangle}

\newcommand		{\intd}			{\int_{\R^d}}
\newcommand		{\iintd}		{\iint_{\R^{2d}}}
\newcommand		{\ka}			{\kappa_{\alpha,\*}}
\renewcommand	{\*}			{{**}}
\newcommand		{\eps}			{\varepsilon}
\newcommand		{\lapfrac}		{\Delta^\frac{\alpha}{2}}
\DeclareMathOperator{\I}		{\mathit{I}}
\newcommand		{\grad}[1]		{\mathrm{D}^{#1}}
\newcommand		{\Dp}[1]		{{\mathfrak{D}}_{#1}}

\usepackage{stmaryrd}
\newcommand		{\subsetArrow}	{\mathrel{\ooalign{$\subset$\cr%
\hidewidth\raise-.087ex\hbox{$_\shortrightarrow\mkern-1.5mu$}\cr}}}
\newcommand		{\subsetarrow}	{\mathrel{\ooalign{$\subset$\cr%
\hidewidth\raise-1.45ex\hbox{$\vec{}\mkern6mu$}\cr}}}

\title{Fractional Fokker-Planck Equation\\ with General Confinement Force}
\author{L. Lafleche}

\def\signll{\bigskip\begin{center}{
	\sc Laurent
	Lafleche\par\vspace{3mm} 
	Université Paris-Dauphine, PSL Research University\par
	CEREMADE, UMR CNRS 7534\par
	Place du Maréchal de Lattre de Tassigny \par
	75775 Paris Cedex 16 FRANCE\par\vspace{3mm} e-mail:}
	\tt{lafleche@ceremade.dauphine.fr}
\end{center}}

\begin{document}

\begin{abstract}This article studies a Fokker-Planck type equation of fractional diffusion with conservative drift
\[
	\partial_tf = \lapfrac f + \divg(Ef),
\]
where $\lapfrac$ denotes the fractional Laplacian and $E$ is a confining force field. The main interest of the present paper is that it applies to a wide variety of force fields, with a few local regularity and a polynomial growth at infinity.

We first prove the existence and uniqueness of a solution in weighted Lebesgue spaces depending on $E$ under the form of a strongly continuous semigroup. We also prove the existence and uniqueness of a stationary state, by using an appropriate splitting of the fractional Laplacian and by proving a weak and strong maximum principle.

We then study the rate of convergence to equilibrium of the solution. The semigroup has a property of regularization in fractional Sobolev spaces, as well as a gain of integrability and positivity which we use to obtain polynomial or exponential convergence to equilibrium in weighted Lebesgue spaces.
\end{abstract}

\maketitle

\begin{center} {\bf 
Version of \today}
\end{center}

\bigskip

Partial differential equations, Asymptotic behavior of solutions
 [see also 45C05, 45K05,  35410],

\bigskip


\renewcommand{\contentsname}{\centerline{Table of Contents}}
\setcounter{tocdepth}{2}	
\tableofcontents


\bigskip
\section{Introduction}
\label{sec:intro}

\subsection{Presentation of the equation and preceding work}

	We consider the homogeneous fractional Fokker-Planck Equation
	\begin{equation*}\label{eq:FFP}\tag{FFP}
		\partial_t f = \Lambda f := \I(f) + \divg\left(E f\right),
	\end{equation*}
	where $E$ is a given force field with polynomial growth at infinity and \begin{equation*}
		\I = \lapfrac \text{ with }  \alpha\in (0,2)
	\end{equation*}
	is the fractional Laplacian. The fractional Laplacian is a generalization of the Laplacian that can be seen as the opposite of a fractional iteration of the positive operator $-\Delta$. It can be defined for any nice function $f$ through its Fourier transform by
	\begin{equation}\label{def:I}
		\widehat{\I(f)} = -|2\pi\xi|^\alpha \widehat{f}.
	\end{equation}
	Alternatively, it is also defined up to a constant depending on $\alpha$ and $d$ for sufficiently smooth functions $f$ by the following integral expression (see e.g. \cite[Chapter 1, \S 1] {landkof_foundations_1972})
	\begin{equation}\label{eq:I_u_intg1}
		\I(f) = \vp\intd \frac{f(y)-f(x)}{|y-x|^{d+\alpha}} \d y,
	\end{equation}
	where $\vp$ indicates that it is a principal value when $\alpha\geq1$.
	
	It can be seen as the infinitesimal generator of a Levy process. A probabilistic point of view about fractional diffusion can for example be found in \cite{jourdain_nonlinear_2008}. The integral representation can be seen in the perspective of the dynamic associated with this Levy process as it represents the fact that particles will jump from $x$ to $y$ proportionally to the difference of value of $f$, from the high to the low densities, and proportionally to the inverse of a power of the distance. It highlights the non-local behavior of this operator.
	
	It is in our case in competition with the force field $E$. For $\alpha<1$, this force field will be stronger in small scales, resulting in possibly discontinuous solutions (see for example \cite{silvestre_regularity_2005}). We restrict ourselves to a force field with at most polynomial growth at infinity.
	
	We mention that another reason for the recent interest about the factional Laplacian is the fact that it can also be seen as a simplified version of the Boltzmann linearized operator, see for example \cite{desvillettes_smoothness_2005}, \cite{mouhot_rate_2006}, \cite{mischler_stability_2009}, \cite{mischler_stability_2009-1}, \cite{tristani_boltzmann_2016}, \cite{canizo_exponential_2016}, \cite{canizo_rate_2017}, \cite{herau_short_2017}. It was for example used extensively in \cite{imbert_weak_2016} and in \cite{silvestre_new_2016} to retrieve Harnack's inequalities and regularity for the Boltzmann equation without cutoff.

\subsection{Main results}\label{ss:MainResults}

	In all this paper, we will denote by $d\in\N^*$ the dimension of the space for the space variable, $\Omega\subset\R^d$ will be an open subset, $\mu$ a measure (or its identification to a Lebesgue measurable function when it is absolutely continuous with respect to the Lebesgue measure) and $m$ a nonnegative weight function which will often be of the form $\weight{x}^k$ for $k\in\R$ where $\weight{x} = \sqrt{1+|x|^2}$. We will often denote by $C$ constants whose exact value have no importance, or write for example $C_a$ when we want to emphasize that the constant depends on $a$, but also use the following notations
	\begin{align*}
		a\lesssim b &\ \overset{\text{def}}{\ssi}\ \exists C>0,\, a\leq Cb
		\\
		a\simeq b &\ \overset{\text{def}}{\ssi}\ a\lesssim b \text{ and } b\lesssim a.
	\end{align*}
	
	Notice that $f,g$ will usually denote functions of time and space while $u,v$ will usually only depend on the space variable $x$. Moreover, $q = p' := \frac{p}{p-1}$ will denote the Hölder conjugate of $p$ and $a\wedge b := \min(a,b)$.
	
	We will mainly work in weighted Lebesgue spaces denoted by $L^p(m)$ for $p \in [1,\infty]$, associated to the norm 
	\begin{equation*}
		\|u\|_{L^p(m)} \ := \ \|um\|_{L^p}.
	\end{equation*}
	We also recall the extension of Sobolev Spaces (see \cite{campanato_proprieta_1963}) to fractional order of derivation, which can be defined through the following semi-norms, generalization of the Hölder property to the Lebesgue spaces for $s\in (0,1)$
	\begin{equation}\label{def:seminorm}
		|u|_{W^{s,p}}^p \ := \ c_{s,d}\iintd \frac{|u(y)-u(x)|^p}{|y-x|^{d+ps}}\d y\d x.
	\end{equation}
	Those are Banach Spaces for the norm $\|u\|^p_{W^{s,p}} := |u|^p_{W^{s,p}} + \|u\|^p_{L^p}$. When $s\in(1,2)$, the norm becomes $\|u\|^p_{W^{s,p}} := \|\nabla u\|^p_{W^{s-1,p}} + \|u\|^p_{L^p}$. See for example \cite{triebel_theory_1992},\cite{triebel_theory_2010},\cite{mazya_sobolev_2011} or \cite{di_nezza_hitchhikers_2012} for a more complete study of these spaces.\vspace{\baselineskip}
	
	We are interested here in a confining force field with polynomial growth taking the form
	\begin{equation}\label{eq:E_ex}
		E \ = \ \langle x\rangle^{\gamma-2} x \ = \ \nabla\left(\frac{\langle x\rangle^{\gamma}}{\gamma}\right),
	\end{equation}
	with $\gamma \in \R$. To simplify the notations, we will sometimes use $\beta := \gamma-2$. The case $E = x = \nabla V(x)$ with $V(x) = \frac{|x|^2}{2}$ is the most studied in the literature (see for example \cite{biler_generalized_2003}, \cite{gentil_levy-fokker-planck_2008}, \cite{gentil_logarithmic_2009}, \cite{tristani_fractional_2015}). In this case the steady state can be computed explicitly and the equation is equivalent up to a scaling to the fractional heat equation (see for example \cite{biler_asymptotics_2001}). Since our method do not use the explicit formula for $E$, we will always assume the following more general hypotheses for a given $\gamma\in\R$.
	
	\paragraph{\textbf{Hypotheses on $E$:}}
	\begin{align}\label{hyp_2bis:grad_E}
		|\nabla E| &\ \lesssim\ \weight{x}^{\gamma-2}
		\\\label{hyp_3:E_confining}
		E\cdot x &\ \gtrsim\ \weight{x}^{\gamma-2}|x|^2.
	\end{align}
	
	Remark also that the kernel in the definition~\eqref{eq:I_u_intg1} of the fractional Laplacian, $\kappa_\alpha : z \mapsto \frac{c_{\alpha,d}}{|z|^{d+\alpha}}$, could be replaced by any symmetric kernel $\kappa_\alpha$ verifying 
	\begin{equation*}
		\kappa_\alpha(z) \simeq \frac{1}{|z|^{d+\alpha}}.
	\end{equation*}
	
	Our first result is about existence and uniqueness of a solution.

	\begin{thm}\label{th:existence}
		Let $m:=\langle x\rangle^k$ with $ k\in(0,\alpha\wedge1)$. Then there exists $p_\gamma > 1$ such that for all $p\in[1,p_\gamma)$, if $f^\mathrm{in}\in L^p(m)$, there exists a unique solution 
		\begin{equation*}
			f\in C^0(\R_+,L^p(m))
		\end{equation*}
		to the \eqref{eq:FFP} equation such that $f(0,\cdot) = f^\mathrm{in}$. Moreover, $\Lambda$ is the generator of a $C^0$-semigroup in $L^p(m)$. 
	\end{thm}
	
	This result generalizes the results obtained by Wei and Tian in \cite{wei_well-posedness_2015}, where the existence was proved for divergence-bounded force fields. The a priori estimates on weighted spaces, from where come the relations between $E$ and $p$, have been already used in the case of the classical Fokker-Planck equation (for example by Gualdani and al in \cite{gualdani_factorization_2013}).
	
	As it can be seen in the proof, to prove the existence of a solution, hypotheses \eqref{hyp_2bis:grad_E} and \eqref{hyp_3:E_confining} can be weakened to the existence of $k\in(0,\alpha\wedge 1)$ and $p>1$ such that
	\begin{align}\nonumber
		E &\ \in\ W^{1,r}_\mathrm{loc} \cap L^\infty_\mathrm{loc} &&\text{for a given } r>2
		\\\nonumber
		E\cdot x &\ \geq\ 0
		\\\label{eq:varphi}
		\varphi_{m,p} &\ := \ \frac{\divg(E)}{q} - E\cdot\frac{\nabla m}{m} \ \leq \ C.
	\end{align}
	In particular, it implies that we do not need to control $|\nabla E|$ but only $\divg(E)$. Moreover, when $\gamma \leq 2$, \eqref{hyp_3:E_confining} is unnecessary.
	
	Remark that when \eqref{hyp_2bis:grad_E} and \eqref{hyp_3:E_confining} hold, then \eqref{eq:varphi} holds for $\gamma\leq2$ or $p$ smaller than a given $p_\gamma\in (1,+\infty)$ which is such that
	\begin{equation}\label{eq:strict_confinement}
		\forall p\in(1,p_\gamma),\,\varphi_{m,p} \ \leq \ b\mathds{1}_\Omega - a \langle x\rangle^{\gamma-2},
	\end{equation}
	for a given $(a,b)\in \R_+^*\times\R$ and a given bounded set $\Omega$. This relation is similar to the Foster-Lyapunov condition for Harris recurrence (see \cite{meyn_stability_1993}, \cite{bakry_rate_2008}, \cite{hairer_asymptotic_2011} and \cite{eberle_quantitative_2016}). When $E$ takes the form \eqref{eq:E_ex}, we can quantify explicitly the value of $p_\gamma = 1 + \frac{k}{d+\gamma-2-k}$.
	
	\begin{thm}\label{th:regu}
		Let $m:=\langle x\rangle^k$ with $ k\in(0,\alpha\wedge1)$ and $f\in L^1(m)$ be a solution to the \eqref{eq:FFP} equation. Then there exists $p_\gamma > 1$ such that $f$ is immediately in all $L^{p}(m)$ for $p<p_\gamma$ and, if $\gamma\leq 2$, $f\in L^\infty(m)$.
	\end{thm}
	
	There has been some recent interest in the regularity theory for integro-differential equations. In \cite{silvestre_holder_2010}, \cite{silvestre_differentiability_2012}, \cite{schwab_regularity_2016}, it is proved that under some regularity conditions on $E$ and if $f\in L^\infty$ is the solution to \eqref{eq:FFP}, then $f$ is actually Hölder continuous or even more differentiable. However, it is also proved in \cite{silvestre_loss_2013} that there can be some loss of regularity when $E$ is not regular enough. As proved in \cite{chamorro_fractional_2016} for divergence free drifts or in Proposition \ref{prop:regu}, we can still obtain fractional Besov or Sobolev regularity in these cases. Theorem \ref{th:regu} gives in particular the regularization from $L^1$ to $L^\infty$ in the case when $E\in C^1_b$, which then allows to use the theorems cited above.

	\begin{thm}\label{th:unicite_equilibre}
		Assume $\gamma>2-\alpha$ and $m=\langle x\rangle^k$ with $0\leq k<\alpha\wedge1$. Then there exists $p^*>1$ such that for any $p\in(1,p^*)$, there exists a unique $F\in L^p(m)\cap L^1_+$ of mass $1$ such that
		\begin{equation*}
			\Lambda F \ = \ 0.
		\end{equation*}
	\end{thm}
	
	This result generalizes the results obtained by Mischler and Mouhot in  \cite{mischler_exponential_2016} and Kavian and Mischler in \cite{kavian_fokker-planck_2015} where it is proved for the classical Laplacian and respectively $\gamma\geq 1$ and $\gamma \leq 1$. It is also close to the result obtained by Mischler and Tristani in \cite{mischler_uniform_2017} where the fractional Laplacian is replaced by integral operators with integrable kernel.
	
	The last and main result is the following rate of convergence towards equilibrium.
	
	\begin{thm}\label{th:cv}
		Assume $\gamma>2-\alpha$ and let $m:=\langle x\rangle^k$ with $0\leq k < (\alpha\wedge1)$. Then, if $\gamma\geq 2$, there exists $a>0$ such that for any $p\in[1,p_\gamma)$,
		\begin{equation*}
			\|f-F\|_{L^p(m)} \ \lesssim \ e^{-at}\|f^\mathrm{in}-F\|_{L^p(m)}.
		\end{equation*}
		If $\gamma\in(2-\alpha,2)$, there exists $p^*>1$ such that for any $p\in(1,p^*)$ and any $\bar{k}<k$, the following rate holds
		\begin{equation*}
			\|f-F\|_{L^p(\bar{m})} \ \lesssim \ \weight{t}^{-\frac{k-\bar{k}}{2-\gamma}}\|f^\mathrm{in}-F\|_{L^p(m)},
		\end{equation*}
		where $\bar{m} = \weight{x}^{\bar{k}}$.
	\end{thm}
	
	This result generalizes the one obtained by Wang in \cite{wang_phi-entropy_2014} where, following the techniques of \cite{gentil_levy-fokker-planck_2008}, exponential convergence of the relative entropy is obtained for force fields $E\in C^1_b$ such that $\forall v\in\R^d, v\cdot\nabla E\cdot v\simeq|v|^2$ and the one obtained by Tristani in \cite{tristani_fractional_2015} where exponential convergence towards equilibrium is proved in $L^p(m)$ in the case $E(x)=x$. It is also the natural extension to the fractional case of the results obtained by Kavian and Mischler in \cite{kavian_fokker-planck_2015} and Mouhot and Mischler in \cite{mischler_exponential_2016}, which correspond respectively to the case $\gamma\in(0,1)$ and $\gamma\geq 1$ for the classical Laplacian. The reason of the lower bound on $\gamma > 2 - \alpha$ is due to the strong nonlocal behavior of the fractional Laplacian which seems to compensate the confining effect of the force field.\vspace{\baselineskip}
	
	The paper is organized as follows. The second section proves some properties of the fractional Laplacian and of the operator $\Lambda$ which will be useful for the various results of the paper.
	
	Section~\ref{sec:existence} proves the existence and uniqueness in the weighted $L^p(m)$ spaces for $p\in(1,2)$. We first create a solution for an approximated problem and then use a priori estimates and compactness properties to obtain a solution to the original problem.
	
	Following the ideas of Nash in \cite{nash_continuity_1958}, section~\ref{sec:gain} of this article generalizes the regularization property of the semigroup associated to the \eqref{eq:FFP} equation as established in \cite{tristani_fractional_2015}. Moreover, a gain of integrability as well as a gain of positivity are also proved, which are useful to deal with convergence without any $L^\infty$ bound.
	
	In section~\ref{sec:steady}, the existence of a stationary state is proved by using an adequate splitting of the operator. It follows the general idea of writing operators as a regularizing part and a dissipative part, as explained in \cite{gualdani_factorization_2013}. We then prove a weak and strong maximum principle and deduce the uniqueness of the equilibrium from the Krein-Rutman Theorem.
	
	The fifth section deals with polynomial convergence when $E$ is not confining enough to create a spectral gap. It uses techniques inspired from \cite{bakry_rate_2008} by using both Foster-Lyapunov estimates introduced by Meyn and Tweedie in \cite{meyn_stability_1993} and a local Poincaré inequality. It proves the first part of Theorem \ref{th:cv}.
	
	Last section is devoted to the proof of the exponential convergence when $E$ is strongly confining (i.e. $\gamma>2$) and follows a different approach as it replaces the use of the Poincaré inequality by the gain of positivity property, following the work of Hairer and Mattingly in \cite{hairer_yet_2011}. It proves the second part of Theorem \ref{th:cv}.

	\paragraph{\bf Acknowledgments:} I wish to acknowledge the help provided by my supervisor, Mr. Stephane Mischler. He gave me very useful advice and I used a lot his course on evolution PDEs \cite{mischler_introduction_2015}. I would also like to thank Ms. Isabelle Tristani and the members of the CEREMADE for their advice.

\section{Main inequalities}

\subsection{Preliminary results about fractional Laplacian}

	We first recall the standard notations that we will use on this paper. We will denote by $\B(E,F)$ the space of continuous linear mappings from $E$ to $F$, by $u_+ := \max(u,0)$ the positive part of $u$. Moreover, we will identify bounded measures on measurable sets of $\R^d$ with bounded radon measures $\mu \in \mathcal{M}(\Omega) := C_0(\Omega)'$ and write
	\begin{align*}
		\int f\mu & \ := \ \int u(x)\mu(\d x), & \mu(A) & \ := \ \int_A \mu,
	\end{align*}
	for any $\mu$-measurable function $u$ and $\mu$-measurable set $A$. We will write the mass of a measure $\langle u\rangle_{\R^d} := \intd u$. We also recall that $\mathcal{D}(\Omega) = C^\infty_c(\Omega)$ and $\mathcal{D}'(\Omega)$ is the space of distributions on $\Omega$. Moreover, we will not write $\Omega$ when $\Omega = \R^d$.
	
	Notice that in order to simplify the computations, we will use the following definition for the power of a vector, $x^a := |x|^{a-1}x$ for any $a\in \R$, and we will use a short notation to simplify the writing of the integrals,
	\begin{align*}
		\ka &:= \kappa_\alpha(x_*-x) & u &:= u(x) & u_* &:= u(x_*),
	\end{align*}  where $x_*$ denote the first variable of integration. We can write for example
	\begin{align*}
		\iint F(u,u_*) = \iint F(u(x),u(x_*)) \d x_* \d x.
	\end{align*}
	
	With these notations and since $\alpha\in(0,2)$, for sufficiently smooth and decaying functions $u$, we can write the fractional Laplacian as a principal value
	\begin{equation*}
		\I(u) \ = \ \vp\left(\intd \ka(u_*-u)\right) \ = \ \lim\limits_{\eps\to0} c_{\alpha,d}\int_{|x-y|>\eps} \frac{u(y)-u(x)}{|y-x|^{d+\alpha}} \d y.
	\end{equation*}
	Remark that the principal value can be removed when $\alpha\in(0,1)$. An other useful expression is 
	\begin{align}
		\I(u) \ &= \ \intd \frac{u(y)-u(x)-(y-x)\nabla u(x)}{|z|^{d+\alpha}} \d z,
		\label{eq:I_u_intg2}
	\end{align}
	By duality, it can also be defined on more general spaces of tempered distributions with a growth smaller than $|x|^{\alpha}$ at infinity by the formula $\langle \I(u),\varphi\rangle_{\mathcal{D}',\mathcal{D}} := \langle u,\I(\varphi)\rangle_{\I(\mathcal{D})',\I(\mathcal{D})}$. In particular, we will mostly use the fractional Laplacian of weight functions of the form $m(x)=\weight{x}^k$ with $k<\alpha$.
	
	Following the model of the Laplacian, we define for $p>1$
	\begin{align}
		\G(u,v) \ &:= \ \intd \frac{\ka}{2}\ (u_*-u)(v_*-v)
		\label{def:prod_grad}\\
		\Dp{p}(u) \ &:= \ \G(u,u^{p-1}) \ \geq \ 0.
		\label{def:Gp}
	\end{align}
	The first quantity can be seen as a generalization of $\nabla u\cdot \nabla v$. It is known as the "Carré du Champs" operator in Probabilities. The second can be seen as a generalization of $\left|\nabla |u|^{p/2}\right|^2$.
	
	The quantity \eqref{def:prod_grad} comes naturally when considering the fractional Laplacian of a product of (sufficiently smooth) functions, since the following formula holds
	\begin{equation}\label{eq:product}
	 	\I(uv) \ = \ u\I(v)+v\I(u)+2\G(u,v).
	\end{equation}
	Moreover, we have the following integration by parts formula
	\begin{equation}\label{eq:ipp}
		\intd u \I(v) \ = \ \intd \I(u) v \ = \ - \intd \G(u,v).
	\end{equation}
	So that in particular, by definition \eqref{def:Gp}
	\begin{equation}\label{eq:I_dissip}
		\intd \I(u)u^{p-1} \ = \ - \intd \Dp{p}(u) \ \leq \ 0.
	\end{equation}
	Remark that these relations also holds when replacing $\kappa_\alpha(x-x_*)$ by a general symmetric kernel $\kappa(x,x_*)$.
	
	It will be useful to remark that the following quantities are equivalent.
	\begin{prop}\label{prop:Gp}
		Let $u$ be such that $\Dp{p}(u)$ is bounded for a given $p\in(1,\infty)$. Then
		\begin{align}\label{eq:Gp_Dp}
			\Dp{p}(u) \ &\simeq  \ \frac{1}{p}\I(|u|^p)-u^{p-1}\I(u)
			\\\label{eq:Gp_Dpp}
			&\simeq  \ \frac{1}{q}\I(|u|^p)-u\I(u^{p-1})
			\\\label{eq:Gp_grad}
			&\simeq  \ \intd \ka |u_*^{p/2}-u^{p/2}|^2,
		\end{align}
		where we recall that $q=p'$ and $a\simeq b$ means here that $a/b$ is bounded by above and below by positive constants depending only on $p$.
	\end{prop}
	
	\begin{demo}[Proposition~\ref{prop:Gp}]
		For the first line, we remark that
		\begin{align*}
			\Dp{p}(u) &= \iintd \ka (u_*-u)(u^{p-1}_* - u^{p-1})
			\\
			&= \iintd \ka\,d_1(u_*/u) |u|^p
			\\
			\frac{1}{p}\I(|u|^p)-u^{p-1}\I(u) &= \frac{1}{p}\iintd \ka (|u_*|^p-|u|^p - p u^{p-1}(u_* - u))
			\\
			&= \frac{1}{p}\iintd \ka\,d_2(u_*/u) |u|^p,
		\end{align*}
		where we recall that $u^p = |u|^{p-1}u$ and we defined for any $z\in\R$,
		\begin{align*}
			d_1(z) &= (z-1)(z^{p-1}-1) \geq 0
			\\
			d_2(z) &= |z|^p - 1 - p (z-1) \geq 0.
		\end{align*}
		Then we remark that $d_1/d_2$ is a bounded positive function since it is continuous on $\R\backslash\{1\}$, converges to $1$ when $|z|\to\infty$ and to $2/p$ when $z\to 1$. Therefore, $d_1\simeq d_2$ and it implies \eqref{eq:Gp_Dp}. The other inequalities are treated in the same way.
	\end{demo}
	
	Another useful result is the estimation of the growth of the fractional Laplacian of weight functions.
	
	\begin{prop}[Fractional Derivation of weight functions]\label{prop_I_m2}
		Let $k\in(0,\alpha\wedge1)$ and $m : x\mapsto\weight{x}^k$ defined for $x\in\R^d$. Then, the following inequality holds
		\begin{equation}\label{eq:I_m2}
			\left|\I(m)\right| \ \leq \ \frac{C}{\weight{x}^{\alpha-k}},
		\end{equation}
		where $C$ is of the form $\frac{C_k\,\omega_d}{(\alpha-k)(2-\alpha)}$. Moreover, when $\alpha<1$ 
		\begin{align}\label{eq:grad_m}
			\grad{\alpha}m \ &\leq \ \frac{C_{\alpha,k}}{\weight{x}^{\alpha-k}},
			\\\label{eq:grad_m2}
			\grad{\alpha}\left(m^{-1}\right) \ &\leq \ \frac{C_{\alpha,|k|}}{\weight{x}^{\alpha}},
		\end{align}
		where $C_{\alpha,k}$ is of the form $\frac{C_k\,\omega_d}{(\alpha-k)(1-\alpha)}$ and $\grad{\alpha}$ is defined by
		\begin{equation}\label{def:grad}
			\grad{\alpha}u \ := \ \intd \ka\,|u_*-u|.
		\end{equation}
	\end{prop}
	
	\begin{demo}[Proposition~\ref{prop_I_m2}]		
		We first look at the case $\alpha\in(0,1)$ and then at the case $\alpha\in(0,2)$ which works only for $\I(m)$.
	
		\step{1. Case $\alpha\in(0,1)$}
		Let $x\in\R^d$ and $R>1$. We split $\grad{\alpha}$ into two parts
		\begin{equation*}
			\grad{\alpha}m \ \leq \ \int_{|x-y|>R} \frac{|m(x)-m(y)|}{|x-y|^{d+\alpha}}\d y + \int_{|x-y|\leq R} \frac{|m(x)-m(y)|}{|x-y|^{d+\alpha}}\d y \ =: \ \mathcal{I}_1 + \mathcal{I}_2.
		\end{equation*}
		For the first part, we remark that since $k\in(0,1)$ and $\forall y\in\R,|\nabla\weight{y}|\leq 1$, we obtain
		\begin{align*}
			|\weight{x}^k-\weight{y}^k| \ \leq \ |\weight{x}-\weight{y}|^k \ \leq \ |x-y|^k.
		\end{align*}
		It leads to
		\begin{equation*}
			\mathcal{I}_1 \ \leq \ \int_{|z|>R} \frac{\d z}{|z|^{d+\alpha-k}}
			\ \leq \ \frac{\omega_d}{(\alpha-k)R^{\alpha-k}}.
		\end{equation*}
		$\bullet$ If $|x|\geq1$, we take $R := |x|/2$. Then $|x|^{-1}\leq \sqrt{2}\weight{x}^{-1}$, from what we deduce
		\begin{equation*}
			\mathcal{I}_1 \ \leq \ \frac{C\,\omega_d}{(\alpha-k)} \frac{1}{\weight{x}^{\alpha-k}}.
		\end{equation*}
		Let $y\in\R^d$ be such that $|x-y|<|x|/2$. For $w\in [x,y] \subset \R^d$, we have $|w| \geq |x| - |x-w| \geq |x|/2$. Thus, we obtain
		\begin{align}
			|m(x)-m(y)| \ & \leq \ |x-y|\sup_{[x,y]}|\nabla m| \label{eq:DL}
			\\
			& \leq \ |x-y|\sup_{w\in[x,y]}|k\weight{w}^{k-2}w|
			\nonumber\\
			& \leq \ 2^{1-k}k\weight{x}^{k-1}|x-y|,
			\nonumber
		\end{align}
		where we used $|x|\leq\weight{x}$ and $\weight{x/2} \geq \weight{x}/2$. It implies the following upper bound
		\begin{equation*}
			\mathcal{I}_2 \ \leq \ C\,\weight{x}^{k-1} \int_{|z|\leq |x|/2} \frac{\d z}{|z|^{d+\alpha-1}} \ \leq \ \frac{C\,\omega_d}{1-\alpha} \weight{x}^{k-\alpha}.
		\end{equation*}
		$\bullet$ If $|x|\leq1$, we take $R := 1$ and we deduce
		\begin{equation*}
			\mathcal{I}_1 \ \leq \ \frac{\omega_d}{(\alpha-k)}.
		\end{equation*}
		Moreover, as $k\weight{x}^{k-1}\leq 1$, \eqref{eq:DL} gives us 
		\begin{equation*}
			|m(x)-m(y)| \ \leq \ |x-y|.
		\end{equation*}
		Therefore
		\begin{equation*}
			\mathcal{I}_2 \ \leq \ \int_{|z|\leq 1} \frac{\d z}{|z|^{d+\alpha-1}} \ \leq \ \frac{\omega_d}{1-\alpha}.
		\end{equation*}
		$\bullet$ We end the proof of \eqref{eq:grad_m} by gathering the two parts together. Since $m\geq 1$, we get $\eqref{eq:grad_m2}$ by remarking that
		\begin{equation*}
			\grad{\alpha}(m^{-1}) \ = \ \intd \ka \left|\frac{m_*-m}{m_*m}\right| \ \leq \ \frac{\grad{\alpha}m}{m}.
		\end{equation*}
		
		\step{2. Proof of \eqref{eq:I_m2}}
		
		We use the integral representation \eqref{eq:I_u_intg2} to change $\mathcal{I}_2$ by
		\begin{equation*}
			\mathcal{I}_2 = \int_{|x-y|\leq R} \frac{|m(x)-m(y)-(x-y)\cdot\nabla m(x)|}{|x-y|^{d+\alpha}}\d y.
		\end{equation*}
		Then \eqref{eq:DL} is replaced by a second order Taylor inequality, which gives
		\begin{equation*}
			|m(x)-m(y)-(x-y)\cdot\nabla m(x)| \ \leq \ C_k\weight{x}^{k-2}|x-y|^2.
		\end{equation*}
		The other parts of the proof are similar to the step $1$.
	\end{demo}

\subsection{Inequalities for the generator of the semigroup}\label{subseq:apriori}

	To get existence, uniqueness and additional gains of weight and regularity on the solutions to the \eqref{eq:FFP} equation, the main inequalities are given in the following
	\begin{prop}\label{prop:estim}
		Let $m=\weight{x}^k$ with $k\in(0,1)$ and $u\in L^p(m\weight{x}^{(\gamma-2)_+})$. If $k<\alpha<1$, the following holds
		\begin{equation}\label{eq:estim_a_priori_1}
			\intd \Lambda(u)u^{p-1}m^p + \intd \Dp{p}(um) \ \leq \ \intd |u|^pm^p\left(\frac{C_k}{\weight{x}^{\alpha-k}}+\varphi_{m,p}\right),
		\end{equation}
		where $\varphi_{m,p}$ is defined by \eqref{eq:varphi} and $\Dp{p}\geq 0$ is defined by \eqref{def:Gp}. If $kp< (\alpha\wedge1)$, we also have
		\begin{equation}\label{eq:estim_a_priori_2}
			\intd \Lambda(u)u^{p-1}m^p + C_p\intd \Dp{p}(um) \ \leq \ \intd |u|^pm^p\left(\frac{C_{k,p}}{\weight{x}^\alpha}+\varphi_{m,p}\right).
		\end{equation}
	\end{prop}
	
	\paragraph{\textbf{Remarks:}} In particular, as already pointed out in introduction, $\varphi_{m,p}$ is always bounded above when $\gamma\leq 2$. When $\gamma>2$, there exists $p_\gamma>1$ such that $\varphi_{m,p}$ is bounded for any $p\in(1,p_\gamma)$. Moreover, in this case, there exists $(a,b)\in\R_+^*\times\R$ such that
	\begin{equation*}
		\varphi_{m,p} \leq b-a\weight{x}^{\gamma-2}.
	\end{equation*}

	Inequality~\eqref{eq:estim_a_priori_2} is more restrictive on $k$ since it needs $k<\alpha/p$, but it has the advantage to work for all $\alpha\in(0,2)$ and to give a second term with a smaller weight.
	
	\begin{lem}\label{lem:Jmp_estimate}
		Let $m=\weight{x}^k$ with $|k|<\alpha\leq 1$ and $u\in L^p(m\weight{x}^{(k-\alpha)/p})$. Then the following inequality holds true
		\begin{eqnarray}\label{eq:Jmp_estimate}
			\left|\intd (\I(mu)-m\I(u))(um)^{p-1}\right| & \leq & C_{k} \left\|u\right\|_{L^p(m\weight{x}^{(k-\alpha)/p})}^p.
		\end{eqnarray}
	\end{lem}
	
	\begin{demo}[Lemma~\ref{lem:Jmp_estimate}]
		Using the integral definition \eqref{eq:I_u_intg1} of $\I$, we have
		\begin{align*}
			\intd j_m(u)v \ & = \ \iintd \ka\, \left((u_*m_*-um)-u(m_*-m)\right)v
			\\
			& = \ \iintd \ka\,\frac{m_*-m}{m_*}(u_*m_*)v.
		\end{align*}
		Thus, by Hölder's inequality, we get
		\begin{equation*}
			\left|\intd j_m(u)v\right| \ \leq \ \left(\iintd \ka\, \frac{|m_*-m|}{m_*}|u_*m_*|^p\right)^\frac{1}{p} \left(\iintd \ka\, \frac{|m_*-m|}{m_*}|v|^q\right)^\frac{1}{q}.
		\end{equation*}
		By the fact that $\frac{|m_*-m|}{m_*} = \frac{|m_*^{-1}-m^{-1}|}{m^{-1}}$ and exchanging $x$ and $x_*$ in the first integral, we obtain
		\begin{equation*}
			\left|\intd j_m(u)v\right| \ \leq \ \left(\intd \frac{\grad{\alpha}m}{m}|um|^p\right)^\frac{1}{p} \left(\intd \frac{\grad{\alpha}(m^{-1})}{m^{-1}}|v|^q\right)^\frac{1}{q}.
		\end{equation*}
		where $\grad{\alpha}$ is defined by \eqref{def:grad}. In particular, if $m=\weight{x}^k$ with $|k|<\alpha$, we obtain from Proposition~\ref{prop_I_m2}
		\begin{equation*}
			\left|\intd j_m(u)v\right| \ \leq \ C_k \left\|u\right\|_{L^p(m\weight{x}^{-\alpha/p})} \left\|v\right\|_{L^q(\weight{x}^{(k-\alpha)/q})},
		\end{equation*}
		which implies \eqref{eq:Jmp_estimate} by taking $v=(um)^{p-1}$.
	\end{demo}

	\begin{demo}[Proposition~\ref{prop:estim}]
		Let $\Phi = \frac{|\cdot|^p}{p}$ and $u\in C_c^\infty$. Then, by definition
		\begin{equation*}
			\intd\Lambda(u)\Phi'(u)m^p\ = \ \intd \I(u)\Phi'(u)m^p+\divg(Eu)\Phi'(u)m^p.
		\end{equation*}
		Let first focus on the term containing the force field $E$.
		We expand the divergence of the product, use the fact that $\Phi'(u)\nabla u = \nabla \Phi(u)$ and integrate by parts the second term to find
		\begin{equation*}
			\intd \divg(E u)\Phi'(u)m^p \ = \ \intd \divg(E)(u\Phi'(u)-\Phi(u))m^p - \Phi(u) E\cdot\nabla m^p.
		\end{equation*}
		By definition of $\Phi$, we obtain
		\begin{equation}\label{eq:estim_a_priori_E}
			\intd \divg(E u)u^{p-1}m^p \ = \ \intd |u|^pm^p \varphi_{m,p},
		\end{equation}
		where $\varphi_{m,p}$ is given by \eqref{eq:varphi}. Let now look at the term containing $\I$. By using \eqref{eq:I_dissip}, we have
		\begin{equation*}
			\intd \I(u)u^{p-1}m^p \ = \ - \intd \Dp{p}(um) + \intd (\I(mu)-m\I(u))(um)^{p-1}.
		\end{equation*}
		By \eqref{eq:Jmp_estimate}, when $|k|<\alpha \leq 1$, we deduce the following inequality for $\I$
		\begin{equation}\label{eq:estim_a_priori_I}
			\intd \I(u)u^{p-1}m^p \ \leq \ - \intd \Dp{p}(um) + C_k \intd |u|^pm^p\weight{x}^{k-\alpha}.
		\end{equation}
		For $\alpha\in (0,2)$, when $kp\in(0,\alpha\wedge1)$, we recall that by relation~\eqref{eq:Gp_Dp},
		\begin{equation*}
			\Dp{p}(u) \ \simeq \ \frac{1}{p}\I(|u|^p)-\I(u)u^{p-1}.
		\end{equation*}
		Hence, using the fractional integration by parts formula \eqref{eq:ipp}, we get
		\begin{equation*}
			\intd \I(u)u^{p-1}m^p \ = \ -C_p\intd \Dp{p}(u)m^p + \frac{1}{p}\intd |u|^p\I(m^p).
		\end{equation*}
		By formula~\eqref{eq:I_m2}, it leads to
		\begin{equation}\label{eq:estim_a_priori_I_2_tmp}
			\intd \I(u)u^{p-1}m^p + C_p\intd \Dp{p}(u)m^p \ \leq \ C_k \intd |u|^pm^p\weight{x}^{-\alpha}.
		\end{equation}
		Now we remark that, by relation~\eqref{eq:Gp_grad}
		\begin{align*}
			\intd \Dp{p}(um) \ &\simeq \ \iintd |(um)_*^{p/2}-(um)^{p/2}|^2
			\\
			&\leq \ 2\iintd |u_*^{p/2}-u^{p/2}|^2m_*^p + |m_*^{p/2}-m^{p/2}|^2|u|^p
			\\
			&\lesssim \ \intd \Dp{p}(u)m^p + \Dp{p}(m)|u|^p.
		\end{align*}
		Moreover, since $\Dp{p}(m) \simeq \G(m^{p/2},m^{p/2})$, by the bound \eqref{eq:ipp}, we obtain
		\begin{align*}
			2\Dp{p}(m) \ \lesssim \I(m^p) + m^{p/2}\left|\I(m^{p/2})\right| \lesssim \frac{m^p}{\weight{x}^\alpha},
		\end{align*}
		where we used \eqref{eq:I_m2} since $kp<\alpha$. Therefore, inequality~\eqref{eq:estim_a_priori_I_2_tmp} becomes
		\begin{equation}\label{eq:estim_a_priori_I_2}
			\intd \I(u)u^{p-1}m^p + C_p\intd \Dp{p}(um) \ \leq \ C_{k,p} \intd |u|^pm^p\weight{x}^{-\alpha}.
		\end{equation}
		We conclude that \eqref{eq:estim_a_priori_1} and \eqref{eq:estim_a_priori_2} hold by combining the inequality for the part with $E$, equation~\eqref{eq:estim_a_priori_E} with the inequalities for the parts with $\I$, \eqref{eq:estim_a_priori_I} and \eqref{eq:estim_a_priori_I_2}.
		All these manipulation can be justified by taking $u_n = \chi_n(\rho_n*u) \to u$ where $\chi_n\in C^\infty_c$ is a cutoff function and $\rho_n\in C^\infty_c$ an approximation of $\delta_0$.
		The main technical point is to obtain an estimate on the following commutator
		\begin{equation*}
			r_n(u) \ := \ (E\cdot\nabla u)*\rho_n - E\cdot\nabla(u*\rho_n).
		\end{equation*}
		In the spirit of DiPerna-Lions commutator estimate (see \cite{diperna_ordinary_1989}) and Lemma~\ref{lem:Jmp_estimate}, we obtain
		\begin{equation*}
			r_n(u) \ \underset{n\to+\infty}{\longrightarrow} 0 \ \mathrm{\ in\ } L^p(m\weight{x}^{-(\gamma-2)_+/p}),
		\end{equation*}
		which ends the proof.
	\end{demo}

\section{Well-posedness}\label{sec:existence}

	This section is devoted to the proof of the part of Theorem~\ref{th:existence} concerning existence and uniqueness of a continuous semigroup. In order to prove the existence of a solution to the \eqref{eq:FFP} equation, we use a viscosity approximation of the equation and a truncation of $E$ and $\I$. We first prove the existence for the approximated problem in $L^2(M)$. We can identify the dual of $V:= H^1(M) = \{u\in L^2(M),\nabla u\in L^2(M)\}$ to $H^{-1}(M)$ by defining $\langle f,g\rangle_{V',V} = \langle fM,gM\rangle_{H^{-1},H^1} = \intd fgM^2$. Moreover, $L^2(M)$ is a Hilbert space for the scalar product $\langle f,g\rangle_{L^2(M)} = \intd fgM^2$. Remark that in the case $\alpha>1$, proving the existence is simpler as the divergence operator is bounded in $H^\alpha$, so that we do not need to use a viscosity approximation.
	
	\begin{lem}[Viscosity Approximation]\label{lem_viscosity}
		Let $M:=\weight{x}^k$ with $k\in\R$ and for $\eps\in(0,1)$ define $\kappa_\alpha^\eps(x) := \kappa_\alpha(x)\mathds{1}_{\{\eps<|x|<1/\eps\}}$ and $\I_\eps(u) := \intd \ka^\eps\left(u_*-u\right)$. Then, there exists a unique solution in 
		\begin{equation*}
			C^0([0,T],L^2(M)) \ \cap \ L^2((0,T),H^1(M)) \ \cap \ H^1((0,T),H^{-1}(M)),
		\end{equation*}
		to the problem
		\begin{equation}
			\partial_tf\ =\ \Lambda_\eps f\ =\ \eps\Delta f + \I_\eps(f) + \divg(E_\eps f),
		\end{equation}
		with $f(0,\cdot)=f^{\mathrm{in}}\in L^2(M)$, $E_\eps\in L^\infty$ and \begin{equation*}
			\left(\divg(E_\eps)-E_\eps\cdot\dfrac{\nabla M^2}{M^2}\right)_+\in L^\infty.
		\end{equation*}
	\end{lem}

	\begin{demo}[Lemma~\ref{lem_viscosity}]
		The result is an application of J.L.Lions Theorem (see for example \cite[Théorème X.9]{brezis_analyse_2005}). We thus prove that the hypotheses of this theorem hold.
		
		\step{1. Continuity of $\Lambda_\eps$}
		
		Let $(f,g)\in H^1(M)^2$. Then
		{\small\begin{align*}
			\langle\Lambda_\eps f,g\rangle_{V',V} & =  \intd-\eps\nabla f \cdot\nabla(g M^2) + \I_\eps(f)g M^2 + \divg(E_\eps f)g M^2
			\\
			& =  \intd\left(-\eps\nabla f\cdot\nabla g - \eps g\nabla f\cdot\dfrac{\nabla M^2}{M^2} + \I_\eps(f)g - E_\eps f\left(\nabla g+g\dfrac{\nabla M^2}{M^2}\right)\right) M^2.
		\end{align*}}Since $\kappa_\alpha^\eps\in L^1$, we can write $I_\eps(f) = \kappa_\alpha^\eps * f - K_\eps f$ where $K_\eps = \|\kappa_\alpha^\eps\|_{L^1}$. Using Peetre's inequality which tells that
		\begin{align*}
			\weight{x+y} \ & \leq \ \sqrt{2}\weight{x}\weight{y},
		\end{align*}
		and the fact that $\kappa_\alpha^\eps$ is compactly supported, we get after a short computation
		\begin{equation}\label{eq:bound_I_eps}
			\left|\intd \I_\eps(f)gM^2\right| 
			\ \leq \ C_\eps K_\eps \|f\|_{L^2(M)}\|g\|_{L^2(M)}.
		\end{equation}
		Thus, using the Cauchy-Schwartz inequality, there exists $C_\eps>0$ such that
		\begin{equation*}
			|\langle\Lambda_\eps f,g\rangle_{V',V}| \ \leq \ \left(C_k(\eps+\|E_\eps\|_{L^\infty})+C_\eps K_\eps\right) \|f\|_{H^1(M)}\|g\|_{H^1(M)},
		\end{equation*}
		where we used $|\nabla M^2| \leq 2|k|M^2$. It proves that $\Lambda_\eps\in\B(V,V')$.
		
		\step{2}
		
		For $f\in H^1(M)$, using \eqref{eq:bound_I_eps} and the a priori estimate \eqref{eq:estim_a_priori_E}, we get
		\begin{align*}
			\langle\Lambda_\eps f,f\rangle_{V',V} & = \intd\left(\I_\eps(f)f -\eps|\nabla f|^2 + f^2 \left(\divg(E_\eps)-E_\eps\cdot\dfrac{\nabla M^2}{M^2}- \eps \dfrac{\Delta M^2}{2M^2}\right)\right) M^2
			\\
			&\leq -\eps\|f\|_{H^1(M)}^2 + C_{k,\eps,E_\eps,\kappa_\alpha^\eps} \|f\|_{L^2(M)}^2,
		\end{align*}
		where we used $|\nabla M^2| \leq 2|k|M^2$ and $|\Delta M^2| \leq 6|k|M^2$. Therefore, we can apply J.L.Lions Theorem.
	\end{demo}
	
	To get results in the good spaces, we will use the following injection that is a straightforward application of Hölder's inequality and the density of $C^\infty_c$ in $L^p$.
	
	\begin{lem}\label{lem_inclusion}
		Let $(p,q)\in[1,+\infty]^2$ and $(l,k)\in\R^2$ such that $p\leq q$ and $(l-k)>d\left(\frac1p-\frac1q\right)$.
		Let $M=\weight{x}^l$ and $m=\weight{x}^k$, then 
		\begin{equation*}
			L^q(M) \hookrightarrow L^p(m),
		\end{equation*}
		with dense and continuous embedding. In particular, if $l>k+\frac d2$ and $p\in[1,2]$, we have the following embedding $L^2(M) \hookrightarrow L^p(m)$.
	\end{lem}

	We now can prove the existence of a weak solution by letting $\eps\to 0$.
	
	\begin{lem}\label{lem:existence}
		Let $m=\weight{x}^k$ with $k\in(0,\alpha\wedge 1)$ and $p\in(1,p_\gamma)$  as defined by \eqref{eq:strict_confinement} (or $p>1$ if $\gamma\leq 2$). Then there exists a unique weak solution $f\in L^\infty_{\mathrm{loc}}(\R_+,L^p(m))$ to the \eqref{eq:FFP} equation.
	\end{lem}
	
	\begin{demo}[Lemma~\ref{lem:existence}]
		We prove first existence of a solution in $L^p(m)$ by using the approximation in $L^2(M)$ and then we use it to prove existence in $L^1(m)$. 
		
		\step{1. Existence in $L^p(m)$ for $p>1$}
		
		Assume that $f^\mathrm{in}\in L^p(m)$ for $p\in(1,2]$. Then, by Lemma~\ref{lem_inclusion}, there exists a family of functions $f_\eps^\mathrm{in}\in L^2(M)$ such that
		\begin{equation*}
			f_\eps^\mathrm{in} \overset{L^p(m)}{\underset{\eps\to 0}{\longrightarrow}} f^\mathrm{in}.
		\end{equation*}
		For a fixed $\eps>0$, let $\chi_\eps\in C^\infty_c$ be a radial function such that $\chi_\eps(x) = \tilde{\chi}_\eps(|x|)$ where $\tilde{\chi}_\eps$ is a decreasing function and $\mathds{1}_{\oball(0,1/\eps)} \leq \chi_\eps \leq \mathds{1}_{\oball(0,2/\eps)}$. Let $f_\eps\in C([0,T],L^2(M))$ be a solution of $\partial_t f_\eps = \Lambda_\eps f_\eps$ as given by Lemma~\ref{lem_viscosity}, with $E_\eps = E\chi_\eps$. For such a $E_\eps$, we have indeed $\divg(E_\eps)-E_\eps\cdot\dfrac{\nabla M^2}{M^2}$ bounded above because of the fact that $E\in L^\infty_\mathrm{loc}$ and $\divg(E)_+\in L^\infty_\mathrm{loc}$.
		
		Let $\rho\in\mathcal{D}(\R^{d+1},\R_+)$ be such that $\int\rho = 1$ and $\mathrm{supp}(\rho)\subset(-1,0)\times\oball(0,1)$ so that  $\rho_n(t,x):=n^{d+1}\rho(nt,n^dx)$ is an approximation of identity. The fractional Laplacian commutes with the convolution by smooth functions (which is an immediate property by using its Fourier definition \eqref{def:I}), thus the regularized function defined by $f_{\eps,n} := f_\eps*\rho_n\in C^\infty(\R_+\times\R^d)\cap L^2(M)$ verifies in the classical sense the equation
		\[
			\partial_t f_{\eps,n} = \Lambda_\eps f_{\eps,n} + r_n,
		\]
		where
		\[
			r_n = (E_\eps\cdot\nabla f_\eps)*\rho_n-E_\eps\cdot\nabla f_{\eps,n}.
		\]
		As proved in \cite[Lemma~II.1]{diperna_ordinary_1989}, since $E_\eps\in  L^1((0,T),W^{1,r}_\mathrm{loc})$ for $r>1$ such that $\frac{1}{p}=\frac{1}{2}+\frac{1}{r}$ and $f_\eps\in L^\infty((0,T),L^2_\mathrm{loc})$, it holds
		\[
			r_n \underset{n\to\infty}{\longrightarrow} 0 \mathrm{\ in \ } L^1((0,T),L^p_{\mathrm{loc}}).
		\]
		Moreover the convergence also holds in $L^1((0,T),L^p(m))$ because $E_\eps$ is compactly supported. Using inequality \eqref{eq:estim_a_priori_1} or \eqref{eq:estim_a_priori_2} for $\I=\I_\eps$ and the fact that $\varphi_{m,p}$ is bounded from above, we obtain
		\begin{equation*}
			\partial_t\left(\intd \frac{|f_{\eps,n}|^p}{p}m^p\right) \ \leq \ \intd |f_{\eps,n}|^pm^p \left(C_k +\frac{\divg(E_\eps)}{q}- k\frac{E_\eps\cdot x}{\weight{x}^2}\right) + |f_{\eps,n}|^{p-1}|r_n|m^p.
		\end{equation*}
		For the part containing $E_\eps$, we have
		\begin{equation*}
			\frac{\divg(E_\eps)}{q} - k\frac{E_\eps\cdot x}{\weight{x}^2} \ = \ \left(\frac{\divg(E)}{q} - k\frac{E\cdot x}{\weight{x}^2}\right)\chi_\eps + \frac{E\cdot\nabla(\chi_\eps)}{q}.
		\end{equation*}
		By hypothesis, the first term is bounded above and the second term is negative since
		\begin{equation}
			E\cdot\nabla(\chi_\eps) \ = \ E\cdot\dfrac{x}{|x|}\tilde{\chi}'(|x|) \ \leq \ 0.
		\end{equation}
		Using Hölder's inequality to control the error term, we obtain
		\begin{align*}
			\partial_t\left(\intd \frac{|f_{\eps,n}|^p}{p}m^p\right) \ &\leq \ C\intd |f_{\eps,n}|^pm^p + \|r_n\|_{L^p(m)} \left(\intd |f_{\eps,n}|^pm^p\right)^{1/p'}
			\\
			&\leq \ \left(C + \|r_n\|_{L^p(m)}\right)\intd |f_{\eps,n}|^pm^p + \|r_n\|_{L^p(m)},
		\end{align*}
		where we used the fact that for $\forall x\geq 0,\ x^{1/p'}\leq 1+x$ since $p'\in[2,\infty)$. Grönwall's inequality gives
		\begin{align*}
			\|f_{\eps,n}\|_{L^p(m)}^p \ &\leq \ e^{CT+p\|r_n\|_{L^1((0,T),L^p(m))}} \left(\|f_{\eps,n}^\mathrm{in}\|_{L^p(m)}^p + p\|r_n\|_{L^1((0,T),L^p(m))}\right).
		\end{align*}
		Passing to the limit in $n$, as $f_{\eps,n}\to f_\eps$ in $L^p(m)$, the error term cancels, hence
		\begin{equation*}
			\|f_\eps\|_{L^p(m)}^p \ \leq \ e^{CT} \|f_\eps^\mathrm{in}\|_{L^p(m)}^p.
		\end{equation*}
		Thus, up to a subsequence, it converges in $\mathcal{D}'([0,T]\times\R^d)$ to $f\in L^\infty([0,T],L^p(m))$.
		Let $\varphi\in\mathcal{D}([0,T]\times\R^d)$. Then $\ka^\eps|\varphi_*-\varphi| \leq \ka|\varphi_*-\varphi|$ which is integrable, thus $\I_\eps(\varphi)$ converges to $\I(\varphi)$ by the Lebesgue dominated convergence Theorem. Therefore, we have
		\begin{equation*}
			\langle f_\eps,\I_\eps(\varphi)\rangle_{\mathcal{D}',\mathcal{D}}
			\ \underset{\eps\to 0}{\longrightarrow} \ \langle f,\I(\varphi)\rangle.
		\end{equation*}
		It implies that $\I_\eps(f_\eps)\underset{\eps\to0}{\longrightarrow} \I(f)$ in $\mathcal{D}'([0,T]\times\R^d)$. We can also easily check that $\divg(E_\eps f_\eps)\underset{\eps\to0}{\longrightarrow} \divg(Ef)$ and $(\partial_t-\eps\Delta)f_\eps \underset{\eps\to0}{\longrightarrow} \partial_tf$ in $\mathcal{D}'(\R_+\times\R^d)$.
		Therefore, we obtain the existence of $f\in L^\infty([0,T],L^p(m))$ verifying the \eqref{eq:FFP} equation. Uniqueness follows directly by remarking that $f^\mathrm{in} = 0 \implies f=0$.
		
		\step{2. Existence in $L^1(m)$}
		
		Consider now the case where $f^\mathrm{in}\in L^1(m)$. As $k<\alpha$, by Lemma~\ref{lem_inclusion} we can find $k<l<\alpha$ and $p\in(1,2)$ such that with $M=\weight{x}^l$, we have $L^p(M)\hookrightarrow L^1(m)$. Let $f_n^\mathrm{in}\underset{n\to\infty}{\longrightarrow} f^\mathrm{in}$ in $L^1(m)$ and $f_n$ be the corresponding solution of the \eqref{eq:FFP} given by the existence in the $L^p$ case. Then, the same proof, but with the $L^1(m)$ estimates, gives
		\begin{eqnarray}\label{eq:cauchy_seq_l1}
			\|f_{n_1}-f_{n_2}\|_{L^1(m)} & \leq & e^{(C_0+C)T} \|f_{n_1}^\mathrm{in}-f_{n_2}^\mathrm{in}\|_{L^1(m)} \underset{n\to\infty}{\longrightarrow} 0.
		\end{eqnarray}
		Therefore, $f_n$ is a Cauchy sequence and we can again verify that it converges to a solution in $L^\infty((0,T),L^1(m))$ of the equation.
	\end{demo}
	
	\begin{lem}\label{lem:w_continuity}
		Let $E\in L^\infty_{\mathrm{loc}}$, $m\in L^0(\R,\R_+^*)$ and $f\in L^\infty((0,T),L^p(m))$ for $p\in(1,+\infty)$ be a weak solution of the \eqref{eq:FFP} equation. Then we have the following continuity in time
		\begin{equation*}
			f\in C^0([0,T],w-L^p(m)),
		\end{equation*}
		where $w-L^p(m)$ indicates that we take the weak topology on $L^p(m)$.
	\end{lem}
	
	\begin{demo}[Lemma~\ref{lem:w_continuity}]
		Let $\varphi\in\mathcal{D}(\R^d)$. As $f$ is solution of \eqref{eq:FFP} in $\mathcal{D}'((0,T)\times\R^d)$, taking $\psi\otimes\varphi \in \mathcal{D}((0,T)\times\R^d)$ as test function, we can write
		\begin{equation*}
			-\int_0^T\intd f(t,x)\partial_t\psi(t)\varphi(x)\d x\d t \ = \ \int_0^T\intd f(t,x)\psi(t)(\I(\varphi)-E\cdot\nabla\varphi)(x)\d x\d t,
		\end{equation*}
		or equivalently
		\begin{equation*}
			\partial_tu_\varphi = v_\varphi\mathrm{\ in\ }\mathcal{D}'(0,T),
		\end{equation*}
		with
		\begin{align*}
			u_\varphi &: t\mapsto\intd f(t,\cdot)\varphi &\mathrm{and}&& v_\varphi &: t\mapsto\intd f(t,\cdot)(\I(\varphi)-E\cdot\nabla\varphi).
		\end{align*}
		For $\varphi\in\mathcal{D}(\R^d)$ and $E\in L^\infty_{\mathrm{loc}}$, we have $\I(\varphi)-E\cdot\nabla\varphi\in L^\infty(\weight{x}^{d+\alpha})$. Thus, as by Lemma~\ref{lem_inclusion}, $f\in L^\infty((0,T),L^p(m)) \subset L^\infty((0,T),L^1(\weight{x}^{-(d+\alpha)}))$, we obtain that $u_\varphi\in L^\infty(0,T)$ and $v_\varphi\in L^\infty(0,T)$. Hence, $u_\varphi\in W^{1,\infty}(0,T)\subset C^0([0,T])$.
		
		Let $p\neq 1$. We now show that the result is still true by replacing $\varphi$ by $g\in L^{p'}(m^{-1})$. First, we remark that $u_g$ is well defined in $L^\infty(0,T)$. Then, by the density of $\mathcal{D}(\R^d)$ in $L^{p'}$, there exists a sequence $(\tilde{\varphi}_n)_{n\in\N}\in\mathcal{D}(\R^d)^\N$ such that $\tilde{\varphi}_n \underset{n\to+\infty}{\longrightarrow}gm^{-1}$ in $L^{p'}$, or equivalently, there exists $\varphi_n := m\tilde{\varphi}_n\in\mathcal{D}(\R^d)$ such that $\varphi_n \underset{n\to+\infty}{\longrightarrow}g$ in $L^{p'}(m^{-1})$. We now look at the sequence of $u_{\varphi_n}$ and write
		\begin{align*}
			\left\|u_{\varphi_n}-u_g\right\|_{C^0([0,T])} \ & = \ \left\|\intd f(t,\cdot)(\varphi_n-g)\right\|_{L^\infty(0,T)}
			\\
			& \leq \ \left\|f\right\|_{L^\infty((0,T),L^p(m))}\left\|\varphi_n-g\right\|_{L^{p'}(m^{-1})} \underset{n\to+\infty}{\longrightarrow} 0.
		\end{align*}
		It proves that $u_g\in C^0(0,T)$.
	\end{demo}
	
	We can now combine the previous lemmas to give the proof of Theorem~\ref{th:existence}.
	
	\begin{demo}[Theorem~\ref{th:existence}]
		Since the time continuity in the weak topology $\sigma(X,X')$ implies the continuity in the strong $X$ topology (see e.g. \cite{engel_one-parameter_1999}), combining Lemmas~\ref{lem:existence} and Lemma~\ref{lem:w_continuity} gives the result in the case $p>1$.
		If $p=1$, we prove the time continuity differently. Using again an $L^p(M)$ approximation sequence $f_n$, we obtain from equation \eqref{eq:cauchy_seq_l1} that it is a Cauchy sequence in $C^0([0,T],L^1(m))$, since
		\begin{align*}
			\|f_{n_1}-f_{n_2}\|_{C^0([0,T],L^1(m))} \ & = \ \sup\limits_{t\in(0,T)} \|f_{n_1}-f_{n_2}\|_{L^1(m)}
			\\
			& \leq \ e^{(C_0+C)T} \|f_{n_1}^\mathrm{in}-f_{n_2}^\mathrm{in}\|_{L^1(m)} \underset{n\to\infty}{\longrightarrow} 0,
		\end{align*}
		from what we conclude that $f\in C^0([0,T],L^1(m))$.
	\end{demo}

\section{Additional properties for solutions to the equation}\label{sec:gain} 

	In this section, we prove that the semigroup associated to the \eqref{eq:FFP} equation actually gives gains of regularity, integrability, weight and positivity, which is useful to retrieve quantitative estimates about the regularity of solutions, to prove uniform in time estimates in weighted Lebesgues spaces and existence and uniqueness of the steady state, as well as quantitative rate of decay towards equilibrium.

\subsection{Gain of regularity and integrability}

	\begin{prop}\label{prop:regu}
		Let $f\in L^1(m)$ be a solution of the \eqref{eq:FFP} equation as given in Theorem~\ref{th:existence} for $m=\weight{x}^k$ with $k\in(0,\alpha\wedge 1)$. Then there exists $c>0$ such that the following inequality holds
		\begin{equation}\label{eq:gain_regu}
			\|f\|_{L^p(m)} \lesssim \left(c+\frac{d(p-1)}{\alpha t}\right)^{\frac{d}{q\alpha}} e^{t\lambda_1} \|f^\mathrm{in}\|_{L^1(m)},
		\end{equation}
		where $\lambda_1$ is the growth bound of $e^{t\Lambda}$ in $L^1(m)$, $q'=p\in[1,p_\gamma)$ and if $\alpha\geq 1$, $p<\alpha/k$.
		Moreover, if $f^\mathrm{in}\in L^p(m)$, we obtain the following Sobolev regularity
		\begin{align}\label{eq:regu_H}
			(fm)^{p/2} \ &\in \ L^2((0,T),H^{\alpha/2}).
		\end{align}
	\end{prop}
	
	\paragraph{\bf Remarks:} Formula \eqref{eq:gain_regu} can also be written in other words
	\begin{eqnarray}
		\|e^{t\Lambda}\|_{L^1(m)\to L^p(m)} & \lesssim & \left(c+t^{\frac{-d}{\alpha q}}\right) e^{t\lambda_1}.
	\end{eqnarray}

	In order to show regularizing properties of the \eqref{eq:FFP} equation, one possibility is to use a fractional variant of the Nash inequality in $L^p(m)$ spaces. In the case of $L^2$ spaces, it is proved for example in \cite[Lemma~5.2]{tristani_fractional_2015}.
	\begin{lem}[Fractional Nash inequality in $L^p(m)$]\label{lem:Nash_frac_m_p}
	Let $p\in[1,2]$ and $m=\weight{x}^k$ with $kp\in(0,\alpha\wedge1)$ or $0<k<\alpha<1$. Then for any $u\in L^p(m)$, we have
		\begin{align}\label{eq:Nash_Lpm_W_alpha}
			\intd \I(u)u^{p-1}m^p \ &\lesssim \ C_{k,p}\left\|u\right\|_{L^p(m)}^p - \left|(um)^\frac{p}{2}\right|_{H^{\alpha/2}}^2
			\\\label{eq:Nash_Lpm_L1}
			&\lesssim \ C_{k,p}\left\|u\right\|_{L^p(m)}^p - \left\|u\right\|_{L^p(m)}^{p+\frac{q\alpha}{d}}  \left\|u\right\|_{L^1(m)}^{\frac{-q\alpha }{d}}.
		\end{align}
	\end{lem}
	
	\begin{demo}[Lemma~\ref{lem:Nash_frac_m_p}]
		By the definition of the Sobolev seminorm \eqref{def:seminorm} and the relation \eqref{eq:Gp_grad}, we remark that
		\begin{equation*}
			\intd \Dp{p}(v) \ \simeq \ |v^\frac{p}{2}|_{H^{\alpha/2}}.
		\end{equation*}
		Therefore, \eqref{eq:Nash_Lpm_W_alpha} is a consequence of inequalities \eqref{eq:estim_a_priori_I} or \eqref{eq:estim_a_priori_I_2}. By using the following Gagliardo-Nirenberg inequalities (see for example \cite{mazya_sobolev_2011})
		\begin{equation*}
				\left\|(um)^{p/2}\right\|_{L^2} \ \lesssim \ \left|(um)^{p/2}\right|_{H^{\frac{\alpha}{2}}}^{\theta} \left\|(um)^{p/2}\right\|_{L^{2/p}}^{1-\theta},
		\end{equation*}
		with $\theta \  = \ \frac{p}{p+q\alpha/d}$, which can also be written
		\begin{equation*}
				\left\|u\right\|_{L^p(m)}^{p/\theta} \ \lesssim \ \left|(um)^{p/2}\right|_{H^{\frac{\alpha}{2}}}^{2} \left\|u\right\|_{L^1(m)}^{p(1/\theta-1)},
		\end{equation*}
		we deduce \eqref{eq:Nash_Lpm_L1} from \eqref{eq:Nash_Lpm_W_alpha}.
	\end{demo}
	
	Nash type inequalities let appear the following family of ordinary differential inequalities that can be solved explicitly and lead to the growth in time given by the following application of Gronwall's inequality.
	
	\begin{lem}\label{lem:ODE_ineq}
		Let $(A,B,C,b)\in \R^4$ and $y\in L^1_+(0,T)$ verifying in the weak sense $\partial_t X \leq B X - A e^{-bCt} X^{1+C}$. Then, the following upper bound holds
		\begin{equation*}
			X \leq \frac{e^{-bt}}{A^{1/C}} \left((B-b)+\frac{1}{Ct}\right)^{1/C}.
		\end{equation*}
	\end{lem}
	
	We can now combine Lemma~\ref{lem:ODE_ineq} with previous Nash type inequalities \eqref{eq:Nash_Lpm_W_alpha} and \eqref{eq:Nash_Lpm_L1} to prove Proposition~\ref{prop:regu}.
	
	\begin{demo}[Proposition~\ref{prop:regu}]
		Let $X = X(t) := \|f\|_{L^p(m)}^p$, $Y := \|f\|_{L^1(m)}^p$ and $\theta := \frac{\alpha}{d(p-1)}>0$. The second fractional Nash inequality~\eqref{eq:Nash_Lpm_L1} can be written
		\begin{equation*}
			\intd \I(f)f^{p-1}m^p \ \leq \ \bar{C}X - \tilde{C} Y^{-\theta}X^{1+\theta}.
		\end{equation*}
		Thus, using the inequality~\eqref{eq:estim_a_priori_E} for the $\divg(E\cdot)$ part of the operator $\Lambda$, we obtain
		\begin{align*}
			\partial_t X \ & = \ p\intd \I(f) f^{q-1} m^p + p\intd f^p m^p \varphi_{m,p}
			\\
			& \leq \ (p\bar{C}+C)X - p\tilde{C} Y^{-\theta}X^{1+\theta}.
		\end{align*}
		Using the fact that $Y\leq e^{q\lambda_1 t}Y(0)$ and Lemma~\ref{lem:ODE_ineq}, we obtain
		\begin{equation*}
			X(t) \ \leq \ e^{q\lambda_1 t} \left(\frac{1}{p\tilde{C}}\right)^{\frac{1}{\theta}} \left(c_p+\frac{1}{\theta t}\right)^{\frac{1}{\theta}} Y(0),
		\end{equation*}
		with $c_p = p\bar{C}+C-q\lambda_1$. It proves \eqref{eq:gain_regu}. Let now $Z := |(fm)^{p/2}|_{H^{\alpha/2}}$ and assume $X(0)$ is bounded. Then by Theorem~\ref{th:existence}, we know that $X\leq e^{tp\lambda_p}X(0)$ for a given $\lambda_p\in\R$. Using now the first fractional Nash inequality~\eqref{eq:Nash_Lpm_W_alpha}, we have
		\begin{equation*}
			\intd \I(f)f^{p-1}m^p \ \leq \ \bar{C}X - Z^q.
		\end{equation*}
		It gives us, by integrating the a priori estimates with respect to time
		\begin{equation*}
			\int_0^TZ^q \ \leq \ X(0) - X(T) + (p\bar{C}+C)\int_0^TX.
		\end{equation*}
		Therefore, we obtain
		\begin{equation*}
			\int_0^T\|(fm)^{p/2}\|_{H^{\alpha/2}}^q \ \leq \ X(0) \left( 1 + (p\bar{C}+C+1)p\lambda_pe^{Tp\lambda_p}\right),
		\end{equation*}
		which gives \eqref{eq:regu_H}.
	\end{demo}
	
\subsection{$L^1(m)\to L^\infty(m)$ Regularization when $\gamma \leq 2$}

	When $\gamma\leq 2$, we have a stronger regularization than Proposition~\ref{prop:regu} since the solutions are globally bounded in space. This property, which will hold also for the equilibrium, will be particularly useful to get the polynomial decay of Theorem~\ref{th:cv}.

	\begin{prop}\label{prop:regu_L_infty}
		Assume $\gamma \leq 2$. Let $f\in L^1(m)$ be a solution of the \eqref{eq:FFP} equation as given in Theorem~\ref{th:existence} with  $m:=\weight{x}^k$ with $2k\in(0,(\alpha\wedge 1))$ or $0\leq k<\alpha\leq1$. Then the following inequality holds
		\begin{equation}\label{eq:gain_regu_L_infty}
			\|f\|_{L^\infty(m)} \ \lesssim \ \left(C+t^{\frac{-d}{\alpha}}\right) e^{\frac{t}{2}(\lambda_1^*+\lambda_1)} \|f^\mathrm{in}\|_{L^1(m)},
		\end{equation}
		where $\lambda_1$ is the growth bound of $e^{t\Lambda}$ in $L^1(m)$, $\lambda_1^*$ the growth bound of $e^{t\Lambda^*}$ and $C\in\R$.
	\end{prop}

	\begin{demo}[Proposition~\ref{prop:regu_L_infty}]
		Since $\gamma \leq 2$, then Theorem~\ref{th:existence} and the inequalities \eqref{eq:estim_a_priori_1} or \eqref{eq:estim_a_priori_2} hold in $L^p(m)$ for all $p\in [1,2]$ and Proposition~\ref{prop:regu} holds for $p=2$. It implies
		\begin{equation}\label{eq:regu_L1_L2}
			\|e^{t\Lambda}\|_{L^1(m)\to L^2(m)} \ \lesssim \ t^{\frac{-d}{2\alpha}}e^{t\lambda_1}.
		\end{equation}
		Moreover, for $g$ solution of the dual equation $\partial_t g = \Lambda^* g := \I(g) - E\cdot\nabla g$, we have
		\begin{align*}
			\intd -(E\cdot \nabla g) g^{p-1}m^{-p} \ & = \ \frac{1}{p}\intd |g|^p\divg(Em^{-p})
			\\
			& = \ \intd |g|^pm^{-p}\left(\frac{\divg(E)}{p} - E\cdot \frac{\nabla m}{m}\right)
			\\
			& \leq \ \frac{\left\|\divg(E)\right\|_{L^\infty}}{p} \intd |g|^pm^{-p},
		\end{align*}
		by combining with formula~\eqref{eq:estim_a_priori_I} that still holds, we obtain the estimate
		\begin{equation*}
			\partial_t\left(\intd |g|^pm^{-p}\right) \ = \ - \intd \Dp{p}(gm^{-1}) + \intd |g|^pm^{-p} \left(C_k+\frac{\left\|\divg(E)\right\|_{L^\infty}}{p}\right).
		\end{equation*}
		Which is the equivalent of \eqref{eq:estim_a_priori_1} for the dual equation in $L^p(m^{-1})$. With the same proof, we get that Theorem~\ref{th:existence} and Proposition~\ref{prop:regu} also hold in $L^p(m^{-1})$ for $p\in [1,2]$, from what we deduce
		\begin{equation}\label{eq:regu_L1_L2_*}
			\|e^{t\Lambda^*}\|_{L^1(m^{-1})\to L^2(m^{-1})} \ \lesssim \ \left(c+t^{\frac{-d}{2\alpha }}\right)e^{t\lambda_1^*},
		\end{equation}
		where $\lambda_1^*$ is the growth bound of $e^{t\Lambda^*}$ in $L^1(m^{-1})$. Since the dual of $L^1(m^{-1})$ and $L^2(m^{-1})$ can be identified with $L^\infty(m)$ and $L^2(m)$, we deduce from \eqref{eq:regu_L1_L2_*} that
		\begin{equation*}
			\|e^{t\Lambda}\|_{L^2(m)\to L^\infty(m)} \ \lesssim \ \left(c+t^{\frac{-d}{2\alpha }}\right)e^{t\lambda_1^*}.
		\end{equation*}
		And combining with \eqref{eq:regu_L1_L2}, by writing $e^{t\Lambda} = e^{\frac{t}{2}\Lambda}e^{\frac{t}{2}\Lambda}$, we end up with
		\begin{equation*}
			\|e^{t\Lambda}\|_{L^1(m)\to L^\infty(m)} \ \lesssim \ \left(C+t^{\frac{-d}{\alpha}}\right) e^{\frac{t}{2}(\lambda_1^*+\lambda_1)},
		\end{equation*}
		which ends the proof.
	\end{demo}

\subsection{Gain of positivity}
	
	We prove in this section the gain and the propagation of strict positivity. It will be useful to prove the uniqueness of the steady state and also, as explained in Proposition \ref{prop:cv_harris}, to get asymptotic estimates when we are not able to prove that the steady state is bounded and use Poincaré inequality. The first proposition is the classical maximum principle.
	
	\begin{prop}[Weak Parabolic Maximum Principle]\label{prop:wPMP} Assume that the conditions of Proposition~\ref{prop:estim} are satisfied and let $f\in L^p(\R_+,L^p(m\weight{x}^{(\gamma-2)_+/p}))$ be such that
		\begin{itemize}
		\item[$\circ$] $(\partial_t -\Lambda)f \geq 0$,
		\item[$\circ$] $f(0,\cdot)=f^\mathrm{in}\geq0$.
		\end{itemize}
	Then $f\geq 0$.
	\end{prop}
	
	\begin{demo}[Proposition \ref{prop:wPMP}]
		Let $g\in L^p(m\weight{x}^{(\gamma-2)_+/p})$, $g_- := (-g)_+$ its negative part and $\Phi(g) := g_+^p$. We remark that
		\begin{equation*}
			\intd \I(g)\Phi'(g)m^p \ \leq \ p\intd \I(g_+)g_+^{p-1}m^p,
		\end{equation*}
		because, as $(g_-)(g_+) = 0$, we have
		\begin{align*}
			-\intd \I(g_-)g_+^{p-1}m^p \ & = \ -\iintd \ka((g_-)_* - g_-)g_+^{p-1}m^p
			\\
			& = \ -\iintd \ka(g_-)_*(g_+^{p-1})m^p \ \leq \ 0.
		\end{align*}
		Thus, if $g$ is such that $\partial_t g \leq \Lambda g$, we get 
		\begin{equation*}
			\partial_t\left(\intd |g_+|^pm^p\right) \ \leq\ p\intd (\Lambda g)g_+^{p-1}m^p \ \leq\ p\intd \Lambda(g_+)g_+^{p-1}m^p.
		\end{equation*}
		Using the a priori estimates \eqref{eq:estim_a_priori_1} or \eqref{eq:estim_a_priori_2}, we obtain
		\begin{equation*}
			\intd |g_+|^pm^p \ \leq \ e^{\lambda t}\intd |g_+^\mathrm{in}|^pm^p.
		\end{equation*}
		We conclude by taking $f=-g$ and remarking that $f^\mathrm{in}_- = 0 \implies f_- =0$.
	\end{demo}

	The second proposition claims that the solutions to the \eqref{eq:FFP} equations are actually bounded by below by a strictly positive function as soon as they have positive mass in a compact set. It implies in particular the strong maximum principle.
	
	\begin{prop}\label{prop:positivity}
		Let $f$ be a solution to the \eqref{eq:FFP} equation with initial condition $f^\mathrm{in}\in L^1_+\cap L^p(m)$. Then for any $a>d+\alpha+\gamma-2$ and $R>0$ sufficiently large, there exists an increasing function $\psi_R\in C^0\cap L^\infty(\R_+^*,\R_+^*)$ such that
		\begin{equation*}
			f(t,x) \ \geq \ \frac{\psi_R(t)}{\weight{x}^a}\int_{B_R} f^\mathrm{in},
		\end{equation*}
		where $B_R$ denotes the ball of size $R$.
	\end{prop}
	
	For a given $r>0$, we define $\chi := \mathds{1}_{B_r}$, $\chi^c := 1-\chi$, $\kappa^c := \kappa_\alpha \chi^c + \kappa_\alpha(r)\chi = \min(\kappa_\alpha,\kappa_\alpha(r))$ and $\kappa := \kappa_\alpha - \kappa^c \geq 0$. As $\kappa^c \in L^1$, we will denote by $K^c := \|\kappa^c\|_{L^1}$ and will decompose $\I$ into
	\begin{align*}
		\I_c(u) \ & := \ \intd \kappa^c_\*\,(u_*-u) \ =\ \kappa^c * u - K^c u
		\\
		\I_\chi(u) \ & := \ \intd \kappa_\*\,(u_*-u) \ =\ \int_{|x-y|<r} \kappa(x-y)(u(y)-u(x))\d y.
	\end{align*}
	Then we define the splitting 
	\begin{equation*}
		\Lambda = A + B,
	\end{equation*}
	where
	\begin{equation*}
		A u = \kappa^c * u \text{ and } B = (I_\chi + \divg(E\ \cdot) - K^c).
	\end{equation*}
	Since the second operator still generates a positive semigroup, the strategy is to use the following Duhamel's formula (see e.g. \cite{arendt_one-parameter_1986})
	\begin{equation*}
		e^{t\Lambda} \ = \ e^{tB} + e^{t\Lambda}\star Ae^{tB},
	\end{equation*}
	where we defined the time convolution of two operators by
	\begin{equation*}
		U\star V : t \mapsto \int_0^t U(t-s)V(s)\d s,
	\end{equation*}
	and to prove that $A$ gives a gain of positivity while $e^{t\Lambda}$ propagates the lower bound. These properties are given in the following lemmas. We will need the following bound by below
	
	\begin{lem}[Bound by below for $\I(m)$]\label{lem:I_m_below}
		Let $m(x):=\weight{x}^k$ with $k<\alpha$.
		Then
		\begin{equation*}
			\I(m) \ \geq \ C_k\weight{x}^{-(d+\alpha)}-\tilde{C}_km.
		\end{equation*}
	\end{lem}
	
	\begin{demo}[Proposition~\ref{lem:I_m_below}]
		We use the above splitting of the fractional Laplacian into $I = I_\chi + I_c$ for $\chi = \mathds{1}_{B_1}$.We first deal with $\I_\chi(m)$ and remark that
		\begin{align}\label{def:lap_tronque}
			\I_\chi(u) \ & := \ \intd \kappa_\*\left(u_*-u-(x_*-x)\cdot\nabla u\right).
		\end{align}	
		By a second order Taylor approximation, for $z\in B_1$, we obtain 
		\begin{align*}
			|m(x+z)-m(x)-z\cdot\nabla m(x)| \ &\leq \ \frac{|z|^2}{2} \left\|\nabla^2 m\right\|_{L^\infty(B_1(x))}.
		\end{align*}
		Thus, by the change of variable $z=x-x_*$ in \eqref{def:lap_tronque}, we can write 
		\begin{align*}
			|\I_\chi(m)| \ & \leq \ \frac{1}{2}\left\|\nabla^2 m\right\|_{L^\infty(B_1(x))} \int_{|z|<1} \frac{\chi(z)\d z}{|z|^{d+\alpha-2}} 
			\\
			& \leq \ \frac{\omega_d}{2(2-\alpha)} \left\|\chi\right\|_{L^\infty}\,\left\|\nabla^2 m\right\|_{L^\infty(B_1(x))}.
		\end{align*}
		In particular, since $m = \weight{x}^k$, we have
		\begin{equation*}
			\left\|\nabla^2 m\right\|_{L^\infty(B_1(x))} \ \leq \ \sup_{|z|<R}|k(|k|+3)\weight{x+z}^{k-2}|.
		\end{equation*}
		Peetre's inequality tells that for all $(x,z)\in\R^{2d}$, we have
		\begin{align*}
			\weight{x+z}^{k-2} \ & \leq \ \sqrt{2}^{|k-2|}\weight{x}^{k-2}\weight{z}^{|k-2|}.
		\end{align*}
		Since $\weight{z}\leq\weight{1}$, we obtain
		\begin{equation}\label{eq:i_loc_m}
			\left|\I_\chi(m)\right| \ \leq \ C_{k} \langle x\rangle^{k-2}.
		\end{equation}
		Now deal with the second part. Since $I_c(m) = \kappa^c*m - K^c m$, we just have to remark that
		\begin{equation*}
			\kappa^c*m\ \geq\ \kappa^c*(m(1)\mathds{1}_{B_1})\ \geq\  \frac{C}{(|x|+1)^{d+\alpha}}.
		\end{equation*}
		Then, by combining with \eqref{eq:i_loc_m}, we obtain
		\begin{equation*}
			\I(m) \ \geq \ \frac{C}{(|x|+1)^{d+\alpha}} - \left( K^c + \frac{C_k}{\weight{x}^2}\right)m.
		\end{equation*}
		what gives the result.
	\end{demo}

	\begin{lem}[Propagation of positivity]\label{lem:propag_positivity}
		Let $f\in L^p((0,T),L^p(m\weight{x}^{{(\gamma-2)_+}/p}))$ be a solution to the \eqref{eq:FFP} equation such that $f^\mathrm{in}>\frac{1}{\weight{x}^a}$ with $a>d+\alpha+\gamma-2$. Then there exists $\lambda>0$ such that
		\begin{eqnarray}\label{eq:propag_positivity}
			f(t,x) & \geq & \frac{e^{-\lambda t}}{\weight{x}^a}.
		\end{eqnarray}
	\end{lem}
	
	\begin{demo}[Lemma~\ref{lem:propag_positivity}]
		Let $\beta := \gamma-2$. We prove that for $\lambda$ large enough, $g(t,x) := \mathfrak{m}(x)\psi(t)$ with $\psi(t) = e^{-\lambda t}$ and $\mathfrak{m}(x) = \weight{x}^k$ with $k<-(d+\alpha+\beta)$ is a subsolution. By Lemma~\ref{lem:I_m_below}, we have, indeed
		\begin{equation*}
			\I(\mathfrak{m}) \ \geq \ \left(C_k\weight{x}^{-(d+\alpha+k)}+\tilde{C}_k\right)\mathfrak{m}.
		\end{equation*}
		We deduce
		\begin{align*}
			(\partial_t - \Lambda)g \ & = \ -\lambda g - \I(\mathfrak{m})\psi(t) - \divg(E\mathfrak{m})\psi(t)
			\\
			& \leq \ \left(-\lambda - C_k\weight{x}^{-(d+\alpha+k)}+\tilde{C}_k - \divg(E) - E\frac{\nabla \mathfrak{m}}{\mathfrak{m}}\right)g
			\\
			& \leq \ \left(\tilde{C}_k-\lambda - C_k\weight{x}^{\beta+\eps} + C\,\weight{x}^\beta\right)g,
		\end{align*}
		where $\eps := -(k+d+\alpha+\beta)>0$. Therefore, by taking $\lambda$ sufficiently large we obtain $(\partial_t - \Lambda)g\leq 0$, i.e. $g$ is a subsolution to the equation. As $g\in L^p_{t,x}(\weight{x}^{\alpha+{\beta_+}/p})$, we can apply the weak parabolic maximum principle, Proposition~\ref{prop:wPMP}, to $f-g$ and we get that $f\geq g$.
	\end{demo}
	
	\begin{lem}[Creation of positivity]\label{lem:lower_bound} For $u\in L^1_+$ the following lower bound holds
		\begin{equation}\label{eq:lower_bound_0}
			\kappa^c * u \ \geq \ \frac{C}{\weight{x}^{d+\alpha}}\int_{B_R} u,
		\end{equation}
		where $C = (\sqrt{2}\,\max(r,R,1))^{-(d+\alpha)}$.
	\end{lem}
	
	\begin{demo}[Lemma~\ref{lem:lower_bound}]
		If $y\in B_R$, then $|x-y|\leq |x|+R$. We deduce the following lower bound
		\begin{equation*}
			\kappa^c *u(x) \ \geq \ \int_{|y|<R} \mathds{1}_{|x-y|<r} \frac{u(y)}{|r|^{d+\alpha}} + \mathds{1}_{|x-y|>r} \frac{u(y)}{(|x|+R)^{d+\alpha}}\d y.
		\end{equation*}
		Let $r_1 := \max(r,R,1)$. As $|x|+R \leq |x|+ r_1$ and $r \leq |x|+ r_1$, we get
		\begin{align*}
			\kappa^c *u(x) \ & \geq \ \frac{1}{(|x|+r_1)^{d+\alpha}} \int_{|y|<R} u(y)\d y \ \geq \ \frac{C}{\weight{x}^{d+\alpha}}\int_{B_R} u,
		\end{align*}
		where $C = (\sqrt{2}\,r_1)^{-(d+\alpha)}$.
	\end{demo}
	
	Now we prove that $e^{tB}$ propagates the fact to have a positive mass in a compact set.
	
	\begin{lem}\label{lem:mass_loc}
		Let $u\in L^1(m)$ and $R>0$. Then for all $\delta > 0$, there exists $\lambda_\delta>0$ such that
		\begin{equation}\label{eq:propag_mass_loc}
			\int_{B_{R+\delta}} e^{tB} u \ \geq \ e^{-\lambda_\delta t}\int_{B_R} u.
		\end{equation}
	\end{lem}
	
	\begin{demo}[Lemma \ref{lem:mass_loc}]
		Let $\eta_0\in C^\infty_c$ be a radially decreasing function such that $\mathds{1}_{B_{\bar{R}}}\leq \eta_0 \leq \mathds{1}_{B_R}$ and $\eta_0>0$ on $B_R$. We also define for all $t>0$, $\eta_t := e^{-\lambda t}\eta_0$ for a given $\lambda > 0$. By construction, this is a subsolution of $\partial_t + E\cdot\nabla$ since
		\begin{equation*}
			\partial_t\eta + E\cdot\nabla\eta \ = \ -\lambda\eta - \left(E\cdot \frac{x}{|x|}\right) |\nabla \eta| \ \leq \  - \lambda \eta.
		\end{equation*}
		Our goal is to prove that for $\lambda$ sufficiently large, we even better have $\partial_t\eta + E\cdot\nabla\eta - I_\chi(\eta) \leq 0$. Therefore, we look at the behaviour of $I_\chi(\eta)$ where $\chi = \mathds{1}_{B_r}$. For $|x|>R$ we have
		\begin{equation*}
			I_\chi(\eta) \ = \ \int_{|x-y|<r} \frac{\eta(y)}{|x-y|^{d+\alpha}} \d y \ \geq \ \frac{1}{r^{d+\alpha}}\int_{B_r(x)} \eta \ \geq \ 0,
		\end{equation*}
		where $B_r(x)$ is the ball of center $x$ and radius $r$.
		In particular, defining $j_R := I_\chi(\eta)(x)$ for $|x|=R$, we have $j_R > 0$. As $\eta\in C^\infty$, we easily deduce $I_\chi(\eta)\in C^\infty$ and the existence of $R' \in (\bar{R},R)$ such that for all $|x|\in[R',R]$, $I_\chi(\eta) \geq j_R/2 > 0$. Therefore, we obtain the following cases
		\begin{align*}
			|x|> R' \ & \implies \ I_\chi(\eta) + \lambda \eta \ \geq \ \lambda \eta \ \geq \ 0
			\\
			|x|< R' \ & \implies \ I_\chi(\eta) + \lambda \eta \ \geq \ \lambda \eta(R') - \|I_\chi(\eta)\|_{L^\infty},
		\end{align*}
		and the latter is positive for $\lambda$ sufficiently large. As $\eta \in C^\infty([0,T]\times B_R)$ all the estimates can easily be made uniform in time and we therefore obtain that
		\begin{equation*}
			(\partial_t-B^*)\eta \ \leq \ 0.
		\end{equation*}
		In particular, by application of the maximum principle (Proposition~\ref{prop:wPMP}) we obtain that $e^{tB^*}\mathds{1}_{B_R} \geq e^{tB^*}\eta_0 \geq \eta \geq e^{-\lambda t}\mathds{1}_{B_{\bar{R}}}$. By the dual definition of positivity, we obtain \eqref{eq:propag_mass_loc}.
	\end{demo}
	
	We can now prove the gain of positivity for the \eqref{eq:FFP} equation.
	\begin{demo}[Proposition~\ref{prop:positivity}]
		We combine \eqref{eq:lower_bound_0} and \eqref{eq:propag_mass_loc} to get
		\begin{equation*}
			Ae^{sB} f^\mathrm{in} \ \geq \ \frac{C_{R,\delta,\chi}e^{-\lambda_\delta s}}{\weight{x}^{d+\alpha}}\,\int_{B_R} f^\mathrm{in},
		\end{equation*}
		where $C = (2\sqrt{2})^{-(d+\alpha)}$. By propagation of the positivity (Lemma~\ref{lem:propag_positivity}), for any $a>d+\alpha+\gamma-2$,
		\begin{equation*}
			e^{(t-s)\Lambda}Ae^{sB} f^\mathrm{in} \ \geq \ \frac{C_{R,\delta,\chi}e^{-\lambda (t-s)}e^{-\lambda_\delta s}}{\weight{x}^a}\,\int_{B_R} f^\mathrm{in}.
		\end{equation*}
		In conclusion, by integrating on $s\in[0,t]$ and using the fact that $e^{tB}\geq 0$ and that by Duhamel's formula
		\begin{equation*}
			e^{t\Lambda} = e^{tB} + e^{t\Lambda}A\star e^{tB} \geq e^{t\Lambda}A\star e^{tB},
		\end{equation*}
		we obtain
		\begin{equation*}
			e^{t\Lambda} f^\mathrm{in} \ \geq \ \frac{\psi(t)}{\weight{x}^a}\,\int_{B_R} f^\mathrm{in}.
		\end{equation*}
		where $\psi(t) = C_{R,\delta,\chi}\frac{e^{-\lambda t} -e^{-\lambda_\delta t}}{\lambda-\lambda_\delta} \in C^0\cap L^\infty(\R_+^*,\R_+^*)$.
	\end{demo}

\section{Existence and uniqueness of the steady state}\label{sec:steady}

\subsection{Splitting of $\Lambda$ as a bounded and a dissipative part}

	This section uses a splitting of the operator $\Lambda$ in a dissipative part in $B\in\B(L^1(m),L^p(m^\theta))$ and a bounded part $A\in\B(L^1,L^1(m))$ in order to bound uniformly in time the solution to the \eqref{eq:FFP} equation and obtain the existence of a steady state. We define the new splitting as $\Lambda = A + B $ with
	\begin{align*}
		A :=M\chi_R \text{ and } B := \Lambda - M\chi_R,
	\end{align*}
	where $M>0$ is a large enough constant and $\mathds{1}_{B_R} \leq \chi_R \leq \mathds{1}_{B_{2R}}$ is a smooth cutoff function.
	
	\begin{prop}\label{prop:estim_b}
		Assume $\beta := \gamma - 2> -\alpha$ and let $k\in\left(0,\alpha\wedge1\right)$ and $p\in(1,p_\gamma)$. Then there exists $\omega:\R_+\to \R_+$ such that
		\begin{equation}\label{eq:estim_b}
			\|e^{tB}\|_{\B(L^p(m),L^p(m^\theta))} \ \lesssim \ \omega(t),
		\end{equation}
		where
		\begin{itemize}
		\item if $\beta \geq 0$, then $\theta=1$ and $\omega(t) = e^{-bt}$,
		\item if $\beta \in (-\alpha,0)$, then $\theta$ is any number in $(0,1]$, $\omega(t) = \weight{t}^{-k(1-\theta)/|\beta|}$ and we require $p<\alpha/k$ if $k>\alpha+\beta$.
		\end{itemize}
		In particular, if $\beta > -\alpha$ and $p_\gamma \geq 1$, there exists $(p,\theta)$ such that $\omega \in L^1(\R_+)$. Moreover, the gain of integrability also holds for $B$ and writes
		\begin{equation}\label{eq:gain_regu_b}
			\|e^{tB}f\|_{\B(L^1(m),L^p(m))} \lesssim t^{-\frac{d}{q\alpha}},
		\end{equation}
		where we recall that $q=p'$.
	\end{prop}
	
	\begin{demo}[Proposition~\ref{prop:estim_b}]
		By inequality~\eqref{eq:estim_a_priori_1}, if $0<kp<\alpha\wedge 1$, we have
		\begin{equation*}
			\frac{1}{p}\partial_t\left(\intd |f|^pm^p\right) \ \leq \ \frac{1}{p}\intd |f|^pm^p \left(\frac{C_k}{\weight{x}^\alpha} + \varphi_{m,p} - M\chi_R\right) - C\intd \Dp{p}(fm).
		\end{equation*}
		Or, by inequality~\eqref{eq:estim_a_priori_2}, we can also get for $k\in(0,\alpha)$,
		\begin{equation*}
			\frac{1}{p}\partial_t\left(\intd |f|^pm^p\right) \ \leq \ \frac{1}{p}\intd |f|^pm^p \left(\frac{C_k}{\weight{x}^{\alpha-k}} + \varphi_{m,p} - M\chi_R\right) - C\intd \Dp{p}(fm).
		\end{equation*}
		From \eqref{eq:strict_confinement}, for $p<p_\gamma$, we have $\varphi_{m,p}\leq b\mathds{1}_\Omega-a\weight{x}^\beta$. Therefore, since $\beta>-\alpha$, if $kp<\alpha$ or $k<\alpha+\beta$, for $M$ and $R$ large enough, we obtain
		\begin{equation}\label{eq:quasi_gap_0}
			\frac{1}{p}\partial_t\left(\intd |f|^pm^p\right) \ \leq \ - a \intd |f|^p \weight{x}^{kp+\beta} - C\intd \Dp{p}(fm),
		\end{equation}
		with $a>0$. In particular
		\begin{equation}\label{eq:B_bound}
			\|e^{tB}\|_{\B(L^p(m))} \ \leq \ 1,
		\end{equation}
		which proves inequality \eqref{eq:estim_b} for small times. If $\beta\geq 0$, since $m^p\leq\weight{x}^{kp+\beta}$, the result immediately follows be Grönwall's inequality. Assume now $\beta <0$ and let $\eps := pk(1-\theta)>0$. By Hölder's inequality, we have
		\begin{equation*}
			\intd |f|^pm^{\theta p} \ \leq \ \left(\intd |f|^p\weight{x}^{\theta kp+\beta}\right)^{\eps/(|\beta|+\eps)}\left(\intd |f|^pm^p\right)^{|\beta|/(|\beta|+\eps)}.
		\end{equation*}
		Combining it with \eqref{eq:quasi_gap_0} (where we replace $k$ by $\theta k$) and \eqref{eq:B_bound} leads to
		\begin{equation}\label{eq:quasi_gap_1}
			\partial_t\left(\intd |f|^p m^{\theta p}\right) \ \leq \ -\bar{a}\left(\intd |f|^pm^{\theta p}\right)^{1+|\beta|/\eps}\left(\intd |f^\mathrm{in}|^pm^p\right)^{-|\beta|/\eps}.
		\end{equation}
		By Grönwall's inequality, we obtain
		\begin{equation*}
			\intd |f|^p m^{\theta p} \ \lesssim \ \frac{1}{t^{\eps/|\beta|}} \intd |f^\mathrm{in}|^pm^p.
		\end{equation*}
		It proves inequality \eqref{eq:estim_b} for large times. Moreover, using this time the second term of the right-hand side in \eqref{eq:quasi_gap_0} and following the same proof as in Proposition~\ref{prop:regu}, we get 
		\begin{equation*}
			\|e^{tB}f\|_{\B(L^1(m),L^p(m))} \lesssim t^{-\frac{d}{q\alpha}} e^{t\lambda_1}.
		\end{equation*}
		But using \eqref{eq:B_bound} for $p=1$ proves that we can take $\lambda_1 = 1$. It concludes the proof.
	\end{demo}
	
\subsection{Existence of a unique steady state}
	
	With this dissipative estimate, the gain of integrability property of Proposition \ref{prop:regu} and the properties of $A$, we obtain the following global in time estimates.
	
	\begin{prop}[Global propagation of $L^p$ norms]\label{prop:asympt}
		Assume $\gamma > 2-\alpha$ and let $f$ be a solution of the \eqref{eq:FFP} under the assumptions of Theorem~\ref{th:existence} with $f^\mathrm{in}\in L^1\cap L^p(m^\theta)$. Then, if $\gamma\geq 2$ and $p<p_\gamma$ or if $\gamma\in(2-\alpha,2)$ and $(p,k)$ is such that there exists $\theta\in (0,1)$ such that Proposition \ref{prop:estim_b} holds with $\omega\in L^1(\R_+)$, there exists $C>0$ such that
		\begin{equation*}
			\|e^{t\Lambda}f^\mathrm{in}\|_{L^p(m^\theta)} \ \leq \ C \left(\|f^\mathrm{in}\|_{L^1} + \|f^\mathrm{in}\|_{L^p(m^\theta)}\right).
		\end{equation*}
	\end{prop}
	
	\begin{demo}[Proposition~\ref{prop:asympt}]
		By noticing that $A\in \B(L^1,L^1(m))$, thanks to Proposition~\ref{prop:estim_b}, we obtain the following sequence of estimates
		\begin{equation*}
			L^1 \underset{e^{t\Lambda}}{\overset{1}{\longrightarrow}} L^1 \underset{A}{\overset{\left\|A\right\|}{\longrightarrow}} L^1(m) \underset{e^{tB/2}}{\overset{\omega_2(t)}{\longrightarrow}} L^p(m) \underset{e^{tB/2}}{\overset{\omega(t)}{\longrightarrow}} L^p(m^\theta),
		\end{equation*}
		where $\omega_2(t) = t^{-d/q\alpha}$ (which is integrable in $0$ since $q > d/\alpha$) and we have indicated the linear operator under the arrow and the corresponding growth rate above the arrows. Hence, by remarking that $\omega\omega_2\in L^1(\R_+)$ using the following Duhamel's Formula
		\begin{equation*}
			e^{t\Lambda} = e^{tB} + e^{tB/2}e^{tB/2} \star Ae^{t\Lambda},
		\end{equation*}
		and the global boundedness of $e^{tB}$ in $L^p(m^\theta)$ given by \eqref{eq:B_bound}, we deduce the announced result.
	\end{demo}

	This proposition together with the positivity properties of the semigroup are sufficient to prove existence and uniqueness of the steady state.
	\begin{demo}[Theorem~\ref{th:unicite_equilibre}]
		Since we have obtained a bound, uniform in time, in the weakly sequentially compact set $L^1_+\cap L^p(m)$, a fixed-point argument allows us to claim the existence of a stationary state. Following the same proof as in \cite[Lemma~3.6]{mischler_linear_2017} or \cite[Theorem 5.1]{kavian_fokker-planck_2015}, we obtain from the previous estimates the existence a stationary state $F\in L^1\cap L^p(m)$ to the \eqref{eq:FFP} equation.
		
		Moreover, by the positivity results obtained in Proposition~\ref{prop:positivity} and since $1\in L^{p'}(m^{-1})\cap L^\infty$ for $p<\dfrac{d}{d-k}$, we obtained the following facts
		\begin{list}{$\bullet$}{}
		\item There exists $F\in L^p(m)\cap L^1_+$ such that $\Lambda F = 0$,
		\item $\Lambda^*1=0$ and $1\in(L^p(m)\cap L^1)'_+$,
		\item $\Lambda$ satisfies the strong and the weak maximum principle.
		\end{list}
		As a consequence of the Krein-Rutman Theorem (see e.g. \cite[Theorem 5.3]{mischler_spectral_2013}), we deduce the uniqueness of a stationary state $F\in L^p(m)\cap L^1_+$ of given mass $\|F\|_{L^1} = \|f^\mathrm{in}\|_{L^1}$. It finishes the proof of Theorem~\ref{th:unicite_equilibre}.
	\end{demo}

\section{Polynomial Convergence to the equilibrium for $\gamma\in (2-\alpha,2)$}

	When $\gamma\in (2-\alpha,2)$, the force field seems not confining enough to get exponential convergence since the derivatives of weighted Lebesgue norms let appear Lebesgue norms with smaller weights. Moreover, when $\gamma < 2-\alpha$, the effect of the force field at infinity is dominated by the effect of the fractional Laplacian, which prevent us from proving any explicit convergence result with our method.

	\subsection{Generalized relative entropy}
	
	In this section, we make a remark about the fact that we can already easily prove a non-quantitative version of the convergence toward equilibrium by generalized entropy method. Assume that there exists a steady state $F>0$ to the $\eqref{eq:FFP}$ equation and let $f$ be a solution of the equation of mass $0$. Then for $h := f/F$, by integration by parts, the following computation formally holds
	\begin{equation*}
		\frac{1}{p}\partial_t\left(\intd |h|^pF\right) \ = \ \intd (\I(h^{p-1})h - h E\cdot \nabla h^{p-1}) F.
	\end{equation*}
	Then, since by formula \eqref{eq:Gp_Dpp}
	\begin{equation*}
		\Dp{p}(h) \ \simeq  \ \frac{1}{q}\I(|h|^p) - h\I(h^{p-1}),
	\end{equation*}
	we get
	\begin{align*}
		\frac{1}{p}\partial_t\left(\intd |h|^pF\right) \ & = \ \frac{1}{q}\intd (\I(|h|^p) - ph^{p-1} E\cdot \nabla h)F - \intd \Dp{p}(h)F
		\\
		& = \ \frac{1}{q}\intd (\I(|h|^p) - E\cdot \nabla |h|^p)F - \intd \Dp{p}(h)F
		\\
		& = \ \frac{1}{q}\intd |h|^p(\I(F)+\divg(EF)) - \intd \Dp{p}(h)F
		\\
		& = \ - \intd \Dp{p}(h)F.
	\end{align*}
	Thus, we obtain
	\begin{eqnarray}\label{eq:dissip_entropy}
		\partial_t\left(\intd |f|^pF^{1-p}\right)
		& \leq & - p\intd \Dp{p}\left(\frac{f}{F}\right)F.
	\end{eqnarray}
	Since $\Dp{p}(h) \geq 0$ and $\Dp{p}(h) = 0 \ssi h$ is constant $\ssi f = F$ (by conservation of the mass), it implies that $\intd |h|^pF$ is a strict Lyapunov functional, which implies the convergence to the equilibrium in $L^p(F^{-1/q})$ (see for example \cite[Chapter 5]{mischler_introduction_2015} or \cite{haraux_systemes_1997}). However, we will prove that with other techniques we will get an explicit rate of convergence.

	\subsection{Fractional Poincaré-Wirtinger inequality}
	
	We prove in this section an inequality looking like a fractional Poincaré-Wirtinger inequality on a bounded set $\Omega$, but for the $p$-dissipation $\Dp{p}$ instead of a fractional gradient and for functions such that the mass is zero on the whole space (i.e. $v$ such that $\langle v\rangle_\mu = 0$).
	
	We define the diameter of $\Omega$ as $\diam(\Omega) := \sup_{(x,y)\in\Omega^2}(|x-y|)$. Moreover, we introduce the following notation for the mass and the $L^p$ norm of a function $u$ for a measure $\mu$,
	\begin{align*}
		\langle u\rangle_{\mu,\Omega} &\ :=\ \frac{1}{\mu(\Omega)}\int_{\Omega} u\mu, & \|u\|_{L^p_\mu(\Omega)}^p &\ :=\ \int_\Omega |u|^p\mu,
	\end{align*}
	and we will use the shortcuts $\langle u\rangle_\mu := \langle u\rangle_{\mu,\R^d}$ and $\|u\|_{L^p_\mu}^p = \|u\|_{L^p_\mu(\R^d)}^p$.
	
	\begin{prop}\label{prop:FPW_loc_p_dissip}
		Let $\mu\in L^\infty_\mathrm{loc}\cap L^1_+$. Then for all $v\in L^p_\mu$ such that $\langle v\rangle_\mu = 0$, for all $\Omega\subset\R^d$ bounded, the following inequality holds
		\begin{align*}
			\int_\Omega \left|v\right|^p\mu \ \leq \ C_\mathrm{PW} \int_\Omega \Dp{p}(v)\mu + \eps_\Omega \left\|v\right\|_{L^p_\mu(\Omega)}^{p-1}\left\|v\right\|_{L^p_\mu(\Omega^c)},
		\end{align*}
		where $C_\mathrm{PW} = \diam(\Omega)^{d+\alpha} \left\|\mu\right\|_{L^\infty(\Omega)}$ and $\eps_\Omega = \frac{\mu(\Omega^c)}{\mu(\Omega)}$.
	\end{prop}
	
	It is a consequence of a the following more natural inequality where we control only the distance to the local mass $\langle u\rangle_{\mu,\Omega}$.
	
	\begin{lem}\label{lem:FPW_loc_p_dissip}
		Let $\Omega\subset\R^d$ be bounded and $\mu\in  L^\infty_+(\Omega)$. Then for all $u\in L^p_\mu$, 
		\begin{equation*}
			0\ \leq\ \int_\Omega u^{p-1}(u-\langle u\rangle_{\mu,\Omega})\mu \ \leq \ C_\mathrm{PW} \int_\Omega \Dp{p}(u)\mu,
		\end{equation*}
		where $C_\mathrm{PW} = \,\mathrm{diam}(\Omega)^{d+\alpha} \frac{\left\|\mu\right\|_{L^\infty(\Omega)}}{\mu(\Omega)}$.
	\end{lem}
	
	\begin{demo}[Lemma~\ref{lem:FPW_loc_p_dissip}]
		We normalize $\mu$ to have $\mu\in\mathcal{P}(\Omega)$ (space of probability measures). For all $u\in L^p_\mu(\Omega)$ the following identity hold
		\begin{align*}
			0\ \leq\ \frac{1}{2}\iint_{\Omega^2} (u_*-u)(u_*^{p-1}-u^{p-1})\mu_*\mu \ & = \ \iint_{\Omega^2} u(u^{p-1}-u_*^{p-1})\mu_*\mu
			\\
			& = \ \iint_{\Omega^2}u^{p-1}(u-u_*)\mu_*\mu
			\\
			& = \ \int_\Omega u^{p-1}(u-\langle u\rangle_{\mu,\Omega})\mu.
		\end{align*}
		Hence, using that $|x-y|<2\,\mathrm{diam}(\Omega)$, we get
		\begin{align*}
			\int_\Omega u^{p-1}(u-\langle u\rangle_{\mu,\Omega})\mu \ & = \ \iint_{\Omega^2} \frac{(u_*-u)(u_*^{p-1}-u^{p-1})}{2|x-x_*|^{d+\alpha}}|x-x_*|^{d+\alpha}\mu_*\mu
			\\
			& \leq \ \mathrm{diam}(\Omega)^{d+\alpha}\|\mu\|_{L^\infty(\Omega)}\intd \Dp{p}(u)\mu.
		\end{align*}
		It concludes the proof.
	\end{demo}
	
	\begin{demo}[Proposition~\ref{prop:FPW_loc_p_dissip}]
		Since $\langle v\rangle_{\mu} = 0$, we have
		\begin{align*}
			\int_\Omega |v|^p\mu \ & = \ \int_\Omega v^{p-1}(v-\langle v\rangle_{\mu,\Omega})\mu + \left(\int_\Omega v^{p-1}\mu\right) \frac{1}{\mu(\Omega)} \left(\int_\Omega v\mu\right)
			\\
			& = \ \int_\Omega v^{p-1}(v-\langle v\rangle_{\mu,\Omega})\mu - \frac{1}{\mu(\Omega)} \left(\int_\Omega v^{p-1}\mu\right)  \left(\int_{\Omega^c} v\mu\right),
		\end{align*}
		and, by using Hölder's inequality, the second term can be bounded in the following way
		\begin{align*}
			\frac{1}{\mu(\Omega)}\left|\left(\int_\Omega v^{p-1}\mu\right)\left(\int_{\Omega^c} v\mu\right)\right| \ & \leq \ \frac{1}{\mu(\Omega)} \left\|v\right\|_{L^p_\mu(\Omega)}^{p-1}\mu(\Omega)^\frac{1}{p}\left\|v\right\|_{L^p_\mu(\Omega^c)}\mu(\Omega^c)^\frac{1}{q}
			\\
			& \leq \ \eps_\Omega \left\|v\right\|_{L^p_\mu(\Omega)}^{p-1}\left\|v\right\|_{L^p_\mu(\Omega^c)}.
		\end{align*}
		We apply Lemma~\ref{lem:FPW_loc_p_dissip} to conclude.
	\end{demo}
	
	\subsection{Lyapunov + Poincaré method}
	
	The following proposition is nothing but the part of Theorem~\ref{th:cv} concerning $\gamma\in(2-\alpha,2)$, leading to polynomial convergence. It is inspired from \cite{bakry_rate_2008} where a Local Poincaré together with a Foster-Lyapunov condition are used in the case of the classical Laplacian to prove convergence in spaces of the form $L^2(F^{-1/2}M)$ where $M$ is an exponential or polynomial weight. As this technique strongly uses the formula for gradient of the product of two functions, which is not available for the fractional Laplacian, we work in spaces of the form $L^p((\lambda F^{1-p} + m^p)^{1/p})$ instead, and we use the fact that $F$ has polynomial decay at infinity.
	
	\begin{prop}\label{prop:cv_poly}
		Assume $\beta := \gamma-2 \in (-\alpha,0)$. Let $m = \weight{x}^k$ and $\bar{m} = \weight{x}^{\bar{k}}$ with $|\beta|<k<\bar{k}<\alpha\wedge1$ and $f\in L^p(\bar{m})$ be a solution to the $\eqref{eq:FFP}$ equation for $p < 1+\frac{k-|\beta|}{d+\alpha-k}$, $p<p_\gamma$ and $p<\frac{\alpha}{k}$ if $k>\alpha-\beta$. Then, the following polynomial convergence holds
		\begin{equation*}
			\|f-F\|_{L^p(m)} \ \lesssim \ \frac{1}{\weight{t}^{(\bar{k}-k)/|\beta|}} \|f^\mathrm{in}-F\|_{L^p(\bar{m})}.
		\end{equation*}
	\end{prop}
	
	\begin{demo}[Proposition~\ref{prop:cv_poly}]
	 	By replacing $f$ by $f-F$ and by conservation of the mass, we can assume $\langle f\rangle_{\R^d} = 0$. By Proposition~\ref{prop:positivity}, $F\gtrsim \frac{c}{\weight{x}^{d+\alpha}}$ and for $p\in(1,p_\alpha)$ with $p_\alpha' = \frac{d+\alpha}{k}$, we have $\eps_0 := kp - (p-1)(d+\alpha) > 0$. Therefore, we have
		\begin{equation}\label{eq:sim}
			F^{1-p} \ \lesssim \ \frac{m^p}{\weight{x}^{\eps_0}},
		\end{equation}
		and we deduce that $f\in L^p(F^{-1/q})$. Moreover $F\in L^1_+$ and $F\in L^\infty(m)$ from Proposition~\ref{prop:regu_L_infty}. Therefore, if $f\in L^p(F^{-1/q})$, by combining the fractional Poincaré-Wirtinger inequality (Proposition~\ref{prop:FPW_loc_p_dissip}) with \eqref{eq:dissip_entropy}, we get for a given $\Omega\subset\R^d$ bounded
		\begin{equation}\label{eq:estim_F}
			C_\mathrm{PW}\,\partial_t\left(\intd |f|^pF^{1-p}\right) \ \leq \ - \int_\Omega |f|^pF^{1-p} + \eps_\Omega \left\|f\right\|_{L^p_{F^{1-p}}(\Omega)}^{p-1}\left\|f\right\|_{L^p_{F^{1-p}}(\Omega^c)}.
		\end{equation}
		Moreover, from estimates \eqref{eq:estim_a_priori_E} and \eqref{eq:estim_a_priori_I_2}, for $kp<\alpha$, we have
		\begin{equation*}
			\frac{1}{p}\partial_t\left(\intd |f|^pm^p\right) \ \leq \ \intd |f|^pm^p \left(\frac{C}{\weight{x}^\alpha} + \varphi_{m,p}\right).
		\end{equation*}
		Or we can also use estimates \eqref{eq:estim_a_priori_E} and \eqref{eq:estim_a_priori_I} to deduce that for $k<\alpha$,
		\begin{equation*}
			\frac{1}{p}\partial_t\left(\intd |f|^pm^p\right) \ \leq \ \intd |f|^pm^p \left(\frac{C}{\weight{x}^{\alpha-k}} + \varphi_{m,p}\right).
		\end{equation*}
		From \eqref{eq:strict_confinement} and one of the two above estimates, if $kp<\alpha$ or $|\beta|<\alpha-k$, we get the Foster-Lyapunov like estimate
		\begin{equation}\label{eq:estim_m}
			\frac{1}{p}\partial_t\left(\intd |f|^pm^p\right) \ \leq \ \intd |f|^p \left(b\mathds{1}_\Omega - a \weight{x}^{kp+\beta}\right),
		\end{equation}
		for a given $(a,b)\in \R_+^2$. We define $M^p := m^p+\lambda C_\mathrm{PW} F^{1-p}$ in order to use the negative part of both estimates \eqref{eq:estim_F} and \eqref{eq:estim_m}. Adding the two expressions, we obtain
		\begin{equation*}
			\partial_t\left(\intd |f|^p M^p\right) \ \leq \ \intd |f|^p \left((b-\lambda F^{1-p})\mathds{1}_{\Omega} + \lambda\eps_\Omega F^{1-p} - a\weight{x}^{kp+\beta}\right).
		\end{equation*}
		Using the fact that $F\in L^\infty(m)$ and \eqref{eq:sim}, we obtain the existence of $c>0$ depending only on $F$ such that
		\begin{equation*}
			\partial_t\left(\intd |f|^p M^p\right) \ \leq \ \intd |f|^p \left((b-\lambda c )\mathds{1}_{\Omega} + \weight{x}^{kp}\left(\lambda\eps_\Omega c\weight{x}^{-\eps_0} - a\weight{x}^{-|\beta|}\right)\right).
		\end{equation*}
		Now we remark that since $|\beta| < k$, for $p < \frac{d+\alpha-|\beta|}{d+\alpha-k} = 1 + \frac{k-|\beta|}{d+\alpha-k}$, we obtain $-\eps_0 < \beta$. Taking also $\lambda > \frac{b}{c}$, we get
		\begin{equation*}
			\partial_t\left(\intd |f|^p M^p\right) \ \leq \ \intd |f|^p\weight{x}^{kp+\beta}\left(\lambda\eps_\Omega c - a\right).
		\end{equation*}
		Thus, by taking $\Omega$ large enough so that $\lambda\eps_\Omega c - a < 0$, we obtain the existence of $\bar{a}>0$ such that
		\begin{equation}\label{eq:quasi_gap}
			\partial_t\left(\intd |f|^p M^p\right) \ \leq \ -\bar{a}\intd |f|^p \weight{x}^{kp+\beta}.
		\end{equation}
		Let $\bar{k} \in (k,\alpha)$ and $\eps := p(\bar{k}-k)>0$. By Hölder's inequality, we have
		\begin{equation*}
			\intd |f|^p\weight{x}^{kp} \ \leq \ \left(\intd |f|^p\weight{x}^{kp+\beta}\right)^{\eps/(|\beta|+\eps)}\left(\intd |f|^p\weight{x}^{\bar{k}p}\right)^{|\beta|/(|\beta|+\eps)}.
		\end{equation*}
		Combining it with \eqref{eq:quasi_gap} and Proposition~\ref{prop:asympt} leads to
		\begin{equation*}
			\partial_t\left(\intd |f|^p M^p\right) \ \lesssim \ -\bar{a}\left(\intd |f|^pM^p\right)^{1+|\beta|/\eps}\left(\intd |f^\mathrm{in}|^p\weight{x}^{\bar{k}p}\right)^{-|\beta|/\eps},
		\end{equation*}
		where we used that by \eqref{eq:sim} and the positivity of $F$, we have $ \weight{x}^{kp} \leq  M^p \lesssim \weight{x}^{kp}$. By Grönwall's inequality, we obtain
		\begin{equation*}
			\intd |f|^p m^p \ \leq \ \intd |f|^p M^p \ \lesssim \ \frac{1}{t^{\eps/|\beta|}} \intd |f^\mathrm{in}|^p\weight{x}^{\bar{k}p},
		\end{equation*}
		which gives the expected result.
	\end{demo}

\section{Exponential Convergence to the equilibrium for $\gamma\geq 2$}
	
	When $\gamma\geq 2$, the confinement is sufficiently strong to get an exponential time decay toward equilibrium for $|x|$ large. To get the local behavior, instead of using a local Poincaré inequality as in previous section, we will use the gain of positivity from Proposition~\ref{prop:positivity}.

	\begin{prop}[Convergence in $L^1(m)$]\label{prop:cv_harris}
		Assume $\gamma \geq 2$ and let $f$ be a solution of the \eqref{eq:FFP} equation with $f^\mathrm{in}\in L^1(m)$. Then, there exists $\bar{a}>0$ such that for any $t\in\R_+$
		\begin{equation*}
			\|f(t)-F\|_{L^1(m)} \ \leq \ e^{-\bar{a}t}\|f^\mathrm{in}-F\|_{L^1(m)}.
		\end{equation*}
	\end{prop}

	\begin{demo}[Proposition~\ref{prop:cv_harris}]
		We want here to use the strategy from Hairer and Mattingly in \cite{hairer_yet_2011} so that we use the following notations $P_t := e^{t\Lambda^*}$, $X := L^1(m)$ and $X' = L^\infty(m^{-1})$ where $m=\weight{x}^k$ with $k\in(0,\alpha\wedge1)$. We recall that from Theorem~\ref{th:existence} we immediately deduce by duality that $P : \R_+ \to \B(X')$ is a positive $C^0$-semigroup such that $P_t 1 = 1$. The strategy consists in proving the following Lyapunov and positivity conditions.	
	
		\step{1. Lyapunov condition} Since $E\cdot x \gtrsim |x|^2$, by Proposition~\ref{prop_I_m2}, we have
		\begin{equation*}
			\Lambda^* m \ = \ \I(m) - E\cdot\nabla m \ \leq \ b - a m.
		\end{equation*}
		Moreover, by using Duhamel's formula, we have
		\begin{equation*}
			e^{(\Lambda^*+a)t} \ = \ m + e^{(\Lambda^*+a)t} \star (\Lambda^*+a) m.
		\end{equation*}
		Therefore, we obtain
		\begin{equation*}
			e^{at}P_t m \ \leq \ m + \int_0^t e^{as}P_s b \d s\ \leq\ m + be^{at}/a,
		\end{equation*}
		from what we deduce
		\begin{equation}\label{eq:markov_lyapunov}
			P_t m \ \leq \ \gamma_tm + c,
		\end{equation}
		with $c = b/a$ and $\gamma_t = e^{-at}\in(0,1)$.
		
		\step{2. Positivity condition} From Proposition~\ref{prop:positivity}, we know that there exists $\nu_t(x) = \nu(t,x)\in L^\infty(\R_+,L^1_+(m))$ strictly positive such that for any $f\in L^1_+(m)$, we have 
		\begin{equation*}
			e^{t\Lambda}f \ \geq \ \nu_t\ \int_{B_R} f.
		\end{equation*}
		By duality, it implies that
		\begin{equation}\label{eq:markov_strict_positivity}
			P_t \geq \langle \nu_t,\cdot\rangle\mathds{1}_{m(x)<r}.
		\end{equation}
		where $r = m(R)$.

		\step{3. Convergence in $L^1(m)$}
	
		We define $m_\lambda := 1+\lambda m$ and the following seminorm on $L^\infty(m^{-1})$
		\begin{equation*}
			|\varphi|_{\dot{L}^\infty(m_\lambda^{-1})} \ := \ \sup_{(x,y)\R^{2d}} \left(\frac{|\varphi(x)-\varphi(y)|}{m_\lambda(x)+m_\lambda(y)}\right).
		\end{equation*}
		Then as proved in \cite[Lemma~2.1]{hairer_yet_2011}, we have
		\begin{equation}\label{eq:seminorm_lambda}
			|\varphi|_{\dot{L}^\infty(m_\lambda^{-1})} \ = \ \inf_{c\in\R}\|\varphi-c\|_{L^\infty(m_\lambda^{-1})}.
		\end{equation}
		Moreover, \cite[Theorem~3.1]{hairer_yet_2011} tells us that since \eqref{eq:markov_lyapunov} and \eqref{eq:markov_strict_positivity} imply that for any fixed time $t>0$ there exists a constant $\bar{\gamma}_t\in(0,1)$ such that
		\begin{equation}\label{eq:markov_contraction}
			|P_t\varphi|_{\dot{L}^\infty(m_\lambda^{-1})} \ \leq \ \bar{\gamma}_t|\varphi|_{\dot{L}^\infty(m_\lambda^{-1})}.
		\end{equation}
		By using the semigroup property, we obtain that the optimal $\bar{a}_t := -\ln(\bar{\gamma}_t) > 0$ verifies $\bar{a}_{t+s}\geq \bar{a}_t+\bar{a}_s$, from what we deduce the existence of $\bar{a}>0$ such that
		\begin{equation*}
			\inf_{c\in\R}\|P_t\varphi-c\|_{X'} \ \leq \ e^{-\bar{a}t}\|\varphi\|_{X'},
		\end{equation*}
		where we replaced $L^\infty(m_\lambda^{-1})$ by $X'$ by equivalence of the norms. Take a sequence $(c_n)_{n\in\N}$ converging to the minimizer. Then, we can write for $f^\mathrm{in}\in L^1(m)$ such that $\langle f^\mathrm{in}\rangle_{\R^d} = 0$
		\begin{align*}
			\langle e^{t\Lambda}f^\mathrm{in},\varphi\rangle_{X,X'} \ & = \ \langle f^\mathrm{in},P_t\varphi-c_n\rangle_{X,X'}
			\\
			& \leq \ \|f^\mathrm{in}\|_{L^1(m)} \|P_t\varphi-c_n\|_{X'}.
		\end{align*}
		Passing to the limit $n\to\infty$, for $f(t) := e^{t\Lambda}f^\mathrm{in}$, we get 
		\begin{equation*}
			\|f\|_{L^1(m)} \ = \ \sup_{\|\varphi\|_{X'}\leq 1}\langle f,\varphi\rangle_{X,X'} \ \leq \ e^{-\bar{a}t}\|f^\mathrm{in}\|_{L^1(m)}.
		\end{equation*}
		Proposition~\ref{prop:cv_harris} follows by taking $f-F$ instead of $f$.
	\end{demo}
	
	\begin{demo}[Theorem~\ref{th:cv}]
		The part concerning polynomial convergence when $\alpha\in(2-\alpha,2)$ was proved in Proposition~\ref{prop:cv_poly}. Therefore we just have to prove the part concerning exponential convergence when $\alpha\geq 2$. Thanks to the regularization property of the semigroup from $L^1(m)$ to $L^p(m)$ as proved in Proposition~\ref{prop:regu}, we know that 
		\begin{equation*}
			\|f-F\|_{L^p(m)} \lesssim \left(c+t^{-\frac{d}{q\alpha}}\right) e^{t\lambda_1} \|f^\mathrm{in}-F\|_{L^1(m)}.
		\end{equation*}
		where $\lambda_1$ is exactly such that
		\begin{equation*}
			\|f-F\|_{L^1(m)} \lesssim e^{t\lambda_1} \|f^\mathrm{in}-F\|_{L^1(m)}.
		\end{equation*}
		From Proposition~\ref{prop:cv_harris}, we deduce that $\lambda_1 = -\bar{a} < 0$, which gives the result.
	\end{demo}


\renewcommand{\bibname}{\centerline{Bibliography}}
\bibliographystyle{abbrv} 
\bibliography{FFP}

\bigskip
\signll

\end{document}